\documentclass[11pt,a4paper,oneside]{article}%
\usepackage{amsfonts}
\usepackage{graphicx}
\usepackage{amsmath,amssymb}
\usepackage{color}
\usepackage{amssymb}
\usepackage{amsmath}
\usepackage[T1]{fontenc}
\usepackage{amsthm}
\usepackage[margin=1in]{geometry}
\usepackage{amssymb}
\usepackage{amsfonts}
\usepackage{graphicx}
\usepackage{amsmath}%
\providecommand{\U}[1]{\protect\rule{.1in}{.1in}}

\newcommand\E{{\mathbb {E}}}

\newtheorem {Lemma}{Lemma}[section]

\newtheorem {Theorem}{Theorem}[section]
\newtheorem {Proposition}{Proposition}[section]
\newtheorem {Corollary}{Corollary}[section]

\newtheorem{Definition}{Definition}[section]

\newtheorem{Remark}{Remark}[section]

\newcommand\F{{\mathcal {F}}}

\newcommand\beq{\begin{equation}}
\newcommand\eeq{\end{equation}}

\def\U{\mathcal{U}}
\def\F{\mathcal{F}}

\begin{document}

\begin{center}
{\Large \bf Moment bounds for dependent sequences in smooth Banach spaces}

\bigskip

J. Dedecker$^a$ and F. Merlev\`{e}de$^b$
\bigskip

\end{center}

$^a$ Universit\'e Paris Descartes, Sorbonne Paris Cit\'e,  Laboratoire MAP5 (UMR 8145).

Email: jerome.dedecker@parisdescartes.fr \bigskip

$^b$ Universit\'{e} Paris Est, LAMA (UMR 8050), UPEM, CNRS, UPEC.

Email: florence.merlevede@u-pem.fr \bigskip

\bigskip


 \bigskip

\begin{center}
\bigskip\textbf{Abstract}
\end{center}
We prove a Marcinkiewicz-Zygmund type inequality
for random variables taking values in a  smooth Banach space.
Next, we obtain some sharp concentration inequalities
for the empirical measure of $\{T, T^2, \cdots, T^n\}$, on a class
of smooth functions, when $T$ belongs to a class of nonuniformly expanding maps
of the unit interval.

\section{Introduction and notations}
Let $({\mathbb B}, |\cdot|_{\mathbb B})$ be a separable Banach space.
The notion of $p$-smooth Banach spaces ($1<p\leq 2$) was introduced in a famous paper
by Pisier (\cite{Pi75}, Section 3). These spaces play the same role with respect to martingales as spaces of type $p$ do with respect to the sums of  independent random variables.

  We shall follow the approach of Pinelis \cite{Pi},
 who showed that
2-smoothness is in some sense equivalent to a control of the second directional derivative of
the map $\psi_2$ defined by
$ \psi_2 (x) = | x |^2_{{\mathbb B}}$. In particular, if
there exists $C>0$ such that, for any $x, u $ in ${\mathbb B}$,
\begin{equation}\label{pinelis}
D^2\psi_2 (x) (u,u)   \leq C |u|^2_{{\mathbb B}}  \, ,
\end{equation}
then the space ${\mathbb B}$ is 2-smooth (here $D^2 g (x)(u,v)$ denotes the second derivative of $g$ at point $x$, in the directions $u, v$). In his 1994 paper, Pinelis
\cite{Pi} used
the property \eqref{pinelis} to derive Burkholder and Rosenthal moment inequlities
as well as exponential
bounds for ${\mathbb B}$-valued martingales.

We shall consider two different classes of $2$-smooth Banach spaces, whose smoothness properties are described as follows.
Let $p$ be a real number in $[2, \infty[$ and let $\psi_p$ be the function from ${\mathbb B}$ to ${\mathbb R}$ defined by
\beq \label{defpsi}
\psi_p (x) = | x |^p_{{\mathbb B}} \, .
\eeq
We say that the separable Banach space $({\mathbb B}, |\cdot |_{\mathbb B})$ belongs to the class ${\ {\mathcal C}}_2 (p, c_p)$ if the function $\psi_p$ is two times  differentiable and satisfies for all $x$ and $u$ in ${\mathbb B}$,
\beq \label{majD2mart}
\big | D^2\psi_p (x) (u,u) \big |  \leq c_p |x|^{p-2}_{{\mathbb B}}|u|^2_{{\mathbb B}}  \, .
\eeq
We say that  $({\mathbb B}, |\cdot |_{\mathbb B})$ belongs to the class ${\widetilde {\mathcal C}}_2 (p, {\tilde c}_p)$ if the  more restricive  inequality holds: for all $x,u,v$ in ${\mathbb B}$,
\beq \label{majD2}
\big | D^2\psi_p (x) (u,v) \big |  \leq {\tilde c}_p |x|^{p-2}_{{\mathbb B}}|u|_{{\mathbb B}} |v|_{{\mathbb B}} \, .
\eeq

Before describing our results, let us quote that the class
${\widetilde {\mathcal C}}_2 (p, {\tilde c}_p)$ contains the ${\mathbb L}^q$-spaces
for $q \geq 2$, for which one can compute the constant $\tilde c_p$. The following
lemma will be proved in Appendix.

\begin{Lemma} \label{constantecp}${}$
\begin{enumerate}
\item For any $q \in [2, \infty[$ and any measure space $({\mathcal X}, {\mathcal A}, \mu)$, the space ${\mathbb L}^q ={\mathbb L}^q ({\mathcal X}, {\mathcal A}, \mu)$ belongs  to the class
${\ {\mathcal C}}_2 (p, c_p)$ with $c_p=p(\max(p,q)-1)$, and to the class
 ${\widetilde {\mathcal C}}_2 (p, {\tilde c}_p)$ with ${\tilde c}_p = p \big ( \max( p, 2q-p) -1\big )$.
\item If ${\mathbb B}$ is a separable Hilbert space then it belongs to the class ${\widetilde {\mathcal C}}_2 (p, {\tilde c}_p)$ with ${\tilde c}_p = p (p-1)$.
\end{enumerate}
\end{Lemma}

The main result of this paper is a Marcinkiewicz-Zygmund type inequality
for the moment of ordrer $p$ of partial sums $S_n$ of ${\mathbb B}$-valued random
variables, when ${\mathbb B}$ belongs to
the class ${\widetilde {\mathcal C}}_2 (p, {\tilde c}_p)$. The upper bound is expressed
in terms of conditional expecations of the random variables with respect to a past $\sigma$-field, and extends
the corresponding upper bound by Dedecker and Doukhan \cite{DD} for real-valued random
variables. As in \cite{Rio00} and \cite{DD}, the proof
is done by writing $\psi_p(S_n)$ as a telescoping sum.  The property \eqref{majD2mart}
enables to use the Taylor integral formula at order 2 to control the terms of
the telescoping sums.

This Marcinkiewicz-Zygmund type bound together with the Rosenthal type bound
given in \cite{DMP} and the deviation inequality given in \cite{DMtpa} provide a complete picture of the moment bounds for sums of ${\mathbb B}$-valued random variables, when ${\mathbb B}$ belongs to
the class ${\widetilde {\mathcal C}}_2 (p, {\tilde c}_p)$.
As we shall see, these bounds apply to a large class of dependent sequences, in the whole range from short to long dependence.

As an application, we shall focus on the ${\mathbb L}^q$-norm of the centered empirical distribution function $G_n$ of the iterates of a nonuniformly expanding map $T$ of the unit interval (modelled by a Young tower with polynomial tails).
On the probability space $[0,1]$ equipped with the $T$-invariant probability $\nu$, the covariance
between two H\"older observables of $T$ and $T^n$ is of order
 $n^{-(1-\gamma)/\gamma}$ for some $\gamma \in (0,1)$. Hence the sequence of the iterates
 $(T^i)_{i \geq 1}$
 is short-range dependent if $\gamma <1/2$ and long-range dependent if
 $\gamma \in [1/2, 1)$.
 Moment and deviation bounds for the ${\mathbb L}^q$-norm of $G_n$ are given in Theorem \ref{appYT1}
 in the short range dependent case, and in Theorems  \ref{appYT2} and \ref{appYT3} in the long range dependent case.  In Remark \ref{remappYT3}, we give some arguments,
 based on a limit theorem for the  ${\mathbb L}^2$-norm of $G_n$,  showing that the deviations bounds of Theorem  \ref{appYT3} are in some sense optimal.

 As a consequence of these results, we obtain in Corollary \ref{appKR}
 a complete picture of the behavior of $\|W_1(\nu_n, \nu)\|_p$ for $p\geq 1$, where
$ W_1(\nu_n, \nu)$ is the Wasserstein distance between the empirical measure
$\nu_n$ of $\{T, T^2, \ldots, T^n\}$  and the invariant distribution $\nu$.
These results are different but complementary to the moment bounds on
$W_1(\nu_n, \nu)-{\mathbb E}(W_1(\nu_n, \nu))$  obtained
by Chazottes and Gou\"ezel \cite{CG} and Gou\"ezel and Melbourne \cite{GoMe}
 as a consequence of a concentration inequality
for separately Lipschitz functionals of $(T, T^2, \ldots, T^n)$. See Section
\ref{Wasser} for a deeper discussion.

All along the paper, the notation
 $a_n \ll b_n$ means that there exists a numerical constant $C$ not
depending on $n$ such that  $a_n \leq  Cb_n$, for all positive integers $n$.

\section{A Marcinkiewicz-Zygmund type inequality}
\setcounter{equation}{0}

Our first result extends Proposition 4 of Dedecker and Doukhan \cite{DD} to smooth Banach spaces belonging to  ${\widetilde {\mathcal C}}_2 (p, {\tilde c}_p)$.

\begin{Theorem} \label{inequality-moment-banach}
Let $p$ be a real number in $[2, \infty[$ and let
$({\mathbb B}, |\cdot |_{\mathbb B})$ be a Banach space belonging to the class ${\widetilde {\mathcal C}}_2 (p, {\tilde c}_p)$. Let $(X_i)_{i \in {\mathbb N}}$ be a sequence of centered random variables in ${\mathbb L}^p ({\mathbb B}) $. Let $({\mathcal F}_i)_{i \geq 0}$ be an increasing sequence of $\sigma$-algebras such that $X_i$ is ${\mathcal F}_i$-measurable, and denote by
$\E_i ( \cdot) = \E ( \cdot | {\mathcal F}_i)$ the conditional expectation with
respect to ${\mathcal F}_i$.
Define then
$$
b_{i,n} = \max_{i \leq \ell \leq n } \Big (\E_0 \Big ( | X_i |^{p/2}_{\mathbb B}  \Big | \sum_{k=i}^{\ell} \E_i ( X_k  ) \Big |^{p/2}_{\mathbb B} \Big ) \Big )^{2/p} \, .
$$
For any integer $n\geq 0$, the following inequality holds:
\beq  \label{IneThinemoment}
\E_0 (| S_n |^p_{\mathbb B} ) \leq K^p \Big ( \sum_{i=1}^n b_{i,n} \Big )^{p/2} \text{almost surely, where $K= \sqrt{2p^{-1} } \sqrt{\max ({\tilde c}_p , p/2)}$}\, .
\eeq
\end{Theorem}
\begin{Remark}
Taking  ${\mathcal F}_0=\{ \Omega , \emptyset \}$,
it follows that, for any integer $n\geq 0$,
\beq  \label{IneThinemomentbis}
\E (| S_n |^p_{\mathbb B} ) \leq K^p \Big ( \sum_{i=1}^n \max_{i \leq \ell \leq n } \Big \Vert | X_i |_{\mathbb B} \Big | \sum_{k=i}^{\ell} \E ( X_k | {\mathcal F}_i ) \Big |_{\mathbb B} \Big \Vert_{p/2} \Big )^{p/2} \text{ where $K= \sqrt{2p^{-1} } \sqrt{\max ({\tilde c}_p , p/2)}$}\, .
\eeq
In addition, if we assume that ${\mathbb P}(|X_k|_{\mathbb B}\leq M)=1$  for
any $k \in \{1, \ldots , n\}$, Inequality \eqref{IneThinemomentbis}  combined with  Proposition \ref{compESnmaxESn} of the appendix leads to the bound
\beq\label{IneThinemomentbismax}
\E \Big(\max_{1 \leq k \leq n}| S_k |^p_{\mathbb B} \Big) \leq C_p M^{p-1}n^{p/2} \Big ( \sum_{k=0}^{n-1} \theta^{2/p}(k) \Big )^{p/2} \, ,
\eeq
where
\[
C_p = \frac 1 2 \Big ( \frac{2p K}{p-1} \Big )^p  + 2^{3p-4} 3^p p \quad \text{and} \quad
\theta(k)= \max \Big \{
 {\mathbb E}(|{\mathbb E}(X_{i}|{\mathcal F}_{i-k})|_{\mathbb B}), i  \in \{k+1, \ldots , n\} \Big \} \, .
\]
A complete proof of Inequality \ref{IneThinemomentbismax} will be given in Section \ref{SectprInmax}.
\end{Remark}
When ${\mathbb B}={\mathbb L}^q$ for $q \geq 2$, the constant $K$ of Inequality \eqref{IneThinemomentbis} is equal to  $\sqrt{\max(4q-2p, 2p) -2}$. However we notice that we can obtain a better constant when the underlying sequence is a martingale differences sequence. More precisely, the following extension of the Marcinkiewicz-Zygmund type inequality obtained by Rio (2009) when the random variables are real-valued holds:
\begin{Theorem}  \label{inequality-moment-banach-martingale}
Let $p$ be a real number in $[2, \infty[$ and let $({\mathbb B}, |\cdot |_{\mathbb B})$ be a Banach space belonging to the class ${ {\mathcal C}}_2 (p, c_p)$. Let $(d_i)_{i \in {\mathbb N}}$ be a sequence of martingale differences with values in ${\mathbb B}$ with respect to an increasing filtration $({\mathcal F}_i)_{i \in {\mathbb N}}$ and such that  for all $i \in {\mathbb N}$, $ \Vert |d_i |_{{\mathbb B}} \Vert_p < \infty$. Then, setting $M_n= \sum_{i=1}^n d_i$, the following inequality holds:
\beq  \label{IneThinemomentbismartgenebanach}
\E (| M_n |^p_{{\mathbb B}} ) \leq (p^{-1} c_p)^{p/2} \Big ( \sum_{i=1}^n \Vert
|d_i|_{{\mathbb B}} \Vert_p^2  \Big )^{p/2} \, .
\eeq
\end{Theorem}
In particular if ${\mathbb B}={\mathbb L}^q ({\mathcal X}, {\mathcal A}, \mu)$ with $q \in [2, \infty[$ and $(T, {\mathcal A}, \nu)$ a measure space, Inequality \eqref{IneThinemomentbismartgenebanach}
combined with Lemma \ref{constantecp} leads to
\beq  \label{IneThinemomentbismart}
\E (| M_n |^p_{q} ) \leq (\max(p,q)-1)^{p/2} \Big ( \sum_{i=1}^n \Vert |d_i |_q \Vert_p^2  \Big )^{p/2} \, ,
\eeq
$| \cdot |_q$ being the norm on ${\mathbb L}^q ({\mathcal X}, {\mathcal A}, \mu)$.

\bigskip

\noindent {\bf Proof of Theorem \ref{inequality-moment-banach}.} As in \cite{Rio00} and \cite{DD}, we shall prove the result by induction. For any $t \in [0,1]$ let
\beq \label{defhnt}
h_n(t) = \E_0 \big (  | S_{n-1} +t X_n |^p_{\mathbb B} \big ) \, .
\eeq
Our induction hypothesis at step $n-1$ is the following: for any $k \leq n-1$,
\beq \label{inductionhypo}
h_k(t) \leq K^p \Big ( \sum_{i=1}^{k-1} b_{i,k}  +t b_{k,k} \Big )^{p/2} \, .
\eeq
Since $K \geq 1$, the above inequality is clearly true for $k=1$. Assuming that it is true for $n-1$, let us prove it at step $n$.

Assume that one can prove that
\beq \label{keyinemoment}
h_n(t) \leq \max( {\tilde c}_p , p/2) \Big (  \sum_{k=1}^{n-1}  b_{k,n}  \int_0^1  \big (  h_k(s) \big )^{1-2/p} ds + b_{n,n} \int_0^t  \big (  h_n(s) \big )^{1-2/p}  ds  \Big )  \, ,
\eeq
then, using our induction hypothesis, it follows that
\begin{multline*}
h_n(t) \leq  \max( {\tilde c}_p , p/2) \Big ( \sum_{k=1}^{n-1} b_{k,n} \int_0^1  K^{p-2} \big (  \sum_{i=1}^{k-1} b_{i,k} + s b_{k,k} \big )^{(p-2)/2} ds + b_{n,n} \int_0^t  \big (  h_n(s) \big )^{1-2/p}  \Big ) \\
\leq  \max( {\tilde c}_p , p/2) \Big ( K^{p-2} \sum_{k=1}^{n-1} b_{k,n} \int_0^1   \big (  \sum_{i=1}^{k-1} b_{i,n} + s b_{k,n} \big )^{(p-2)/2} ds + b_{n,n} \int_0^t  \big (  h_n(s) \big )^{1-2/p} ds \Big ) \, .
\end{multline*}
Integrating with respect to $s$, we get
$$
b_{k,n} \int_0^1   \Big (  \sum_{i=1}^{k-1} b_{i,n} + s b_{k,n} \Big )^{(p-2)/2} ds = \frac{2}{p} \Big (  \sum_{i=1}^{k} b_{i,n} \big )^{p/2} - \frac{2}{p} \Big (  \sum_{i=1}^{k-1} b_{i,n} \Big )^{p/2} \, ,
$$
implying that
$$
\sum_{k=1}^{n-1} b_{k,n} \int_0^1   \Big (  \sum_{i=1}^{k-1} b_{i,n} + s b_{k,n} \Big )^{(p-2)/2} ds = 2p^{-1} \Big (  \sum_{i=1}^{n-1} b_{i,n} \Big )^{p/2} \, .
$$
Therefore, since $K^2= 2p^{-1} \max ({\tilde c}_p , p/2)$,
\beq \label{inequation1}
h_n(t)  \leq  K^{p} \Big (  \sum_{i=1}^{n-1} b_{i,n} \Big )^{p/2} +  \max( {\tilde c}_p , p/2) b_{n,n} \int_0^t  \big (  h_n(s) \big )^{1-2/p} ds  \, .
\eeq
Let $H_n (t) = \int_0^t  \big (  h_n(s) \big )^{1-2/p} ds $. The differential integral inequation \eqref{inequation1} writes
$$
H'_n(s) \Big (  K^{p} \Big (  \sum_{i=1}^{n-1} b_{i,n} \Big )^{p/2} +  \max( {\tilde c}_p , p/2) b_{n,n}  H(s) \Big )^{-1+2/p} \leq 1 \, .
$$
Setting $$
R_n(s) = \Big (  K^{p}\Big (  \sum_{i=1}^{n-1} b_{i,n} \Big )^{p/2} +  \max( {\tilde c}_p , p/2) b_{n,n}  H(s) \Big )^{2/p}\, ,
$$
the previous inequality can be rewritten as
$$
R_n'(s) \leq 2 p^{-1}\max( {\tilde c}_p , p/2) b_{n,n}  \, .
$$
Integrating between $0$ and $t$, we derive
\[
\big ( h_n(t) \big )^{2/p} -  K^{2}  \sum_{i=1}^{n-1} b_{i,n}  \leq  R_n(t) - R_n(0)  \leq 2  t p^{-1}  \max( {\tilde c}_p , p/2) b_{n,n} \, .
\]
Taking into account that $K^2= 2p^{-1} \max ({\tilde c}_p , p/2)$, it follows that
\[
\big ( h_n(t) \big )^{2/p} \leq   K^{2}  \Big (  \sum_{i=1}^{n-1} b_{i,n} + t b_{n,n}
\Big ) \, ,
\]
showing that our induction hypothesis holds true at step $n$. To end the proof it suffices to prove \eqref{keyinemoment}. We shall proceed as in the proof of Theorem 2.3 in
\cite{Rio00}. With this aim, let
\[
S_n(t) = \sum_{i=1}^n Y_i(t)\, , \quad
\text{ where $Y_i (t) =X_i$ for $1 \leq i \leq n-1$ and $Y_n(t) = tX_n$.}
\]
Notice that for any integer $k$ in $[1,n-1]$, $S_{k}(t) = S_k$.
Let now $\psi_p$ be defined by \eqref{defpsi}. Applying Taylor integral formula
at order 2, we get
\begin{multline*}
\psi_p( S_n(t)) = \sum_{i=1}^n \big ( \psi_p( S_i(t)) - \psi_p (S_{i-1}(t)) \big ) \\
 = \sum_{k=1}^n D \psi_p( S_{k-1})  (Y_k (t)) + \sum_{i=1}^n \int_0^1 (1-s) D^2 \psi_p (S_{i-1} + s Y_i(t)) (Y_i(t), Y_i(t) ) ds \, .
\end{multline*}
But, for any integer $k$ in $[1,n]$,
\begin{multline*}
D \psi_p( S_{k-1})  (Y_k (t))  = \sum_{i=1}^{k-1} \big ( D \psi_p( S_{i})  (Y_{k} (t))-D \psi_p( S_{i-1})  (Y_k (t))\big ) \\
= \sum_{i=1}^{k-1} \int_0^1 D^2 \psi_p (S_{i-1} + s X_i) (Y_k(t), X_i ) ds \, .
\end{multline*}
Notice now that for any $x$ and $u$ in ${\mathbb B}$, $D^2 \psi_p (x) (u, u) \geq 0$. Indeed, the function $x \mapsto \psi_p(x)= |x|_{\mathbb B}^{p/2}$ is convex for any $p \geq 2$ and is by assumption $2$-times  differentiable, implying  that the second differentiable derivative at $x$ in the direction $u$ is non-negative. So, overall, using the fact that $D^2 \psi_p (x) (u, u) \geq 0$,
\begin{multline*}
\psi_p( S_n(t)) \leq \sum_{i=1}^{n-1} \int_0^1 D^2 \psi_p (S_{i-1} + s X_i) \Big( \sum_{k=i+1}^n Y_k(t), X_i \Big) ds
\\
+  \sum_{i=1}^n \int_0^1 D^2 \psi_p (S_{i-1} + s Y_i(t)) (Y_i(t), Y_i(t) ) ds \, .
\end{multline*}
Taking the conditional expectation w.r.t. ${\mathcal F}_0$ and recalling the definition \eqref{defhnt} of $h_n(t)$, it follows that, for any $t \in [0,1]$,
\begin{multline*}
h_n(t) \leq \sum_{i=1}^{n-1} \int_0^1 \E_0 \Big (D^2 \psi_p (S_{i-1} + s X_i) \Big( \sum_{k=i}^{n-1} X_k +t X_n, X_i \Big) ds \Big ) \\
+ t^2 \int_0^1 \E_0 \Big ( D^2 \psi_p (S_{n-1} + s t X_n) (X_n,X_n ) ds \Big ) \, .
\end{multline*}
Using again the fact that $D^2 \psi_p (v) (u, u) \geq 0$, we have
$$
 t^2 \int_0^1 \E_0 \Big ( D^2 \psi_p (S_{n-1} + s t X_n) (X_n,X_n ) ds \Big ) \leq   \int_0^t \E_0 \Big ( D^2 \psi_p (S_{n-1} + u X_n) (X_n,X_n ) du \Big ) \, .
$$
Hence setting
$$
a_{i,n} (t) = X_i+\sum_{k=i+1}^{n-1} \E (X_k | {\mathcal F}_i) + t \E (X_n | {\mathcal F}_i) \, ,
$$
and using the fact that $({\mathcal F}_i)$ is an increasing sequence of $\sigma$-algebras, we derive
\[
h_n(t) \leq \sum_{i=1}^{n-1} \int_0^1 \E_0 \Big (D^2 \psi_p (S_{i-1} + s X_i) ( a_{i,n} (t), X_i ) ds \Big )
+ \int_0^t \E_0 \Big ( D^2 \psi_p (S_{n-1} + s  X_n) (X_n,X_n ) ds \Big ) \, .
\]
Using \eqref{majD2}, we then get
\begin{multline*}
h_n(t) \leq {\tilde c}_p \sum_{i=1}^{n-1} \int_0^1 \E_0 \Big (|S_{i-1} + s X_i|^{p-2}_{\mathbb B} |  a_{i,n} (t)|_{\mathbb B}| X_i|_{\mathbb B}  \Big ) ds
+ {\tilde c}_p \int_0^t \E_0 \Big (|S_{n-1} + s X_n|^{p-2}_{\mathbb B} | X_n|^2_{\mathbb B}  \Big ) ds \, .
\end{multline*}
H\"older's inequality implies that
\begin{multline} \label{avantdernier}
h_n(t) \leq {\tilde c}_p \sum_{i=1}^{n-1} \int_0^1 \big ( h_{i}(s) \big )^{(p-2)/p} \Big ( \E_0 \big (  |  a_{i,n} (t)|^{p/2}_{\mathbb B}| X_i|^{p/2}_{\mathbb B} \big ) \Big )^{2/p} ds \\
+ {\tilde c}_p \int_0^t \big ( h_{n}(s) \big )^{(p-2)/p} \Big ( \E_0 (| X_n|^p_{\mathbb B} ) \Big )^{2/p} ds \, .
\end{multline}
Let $G_{i,n}(t) = \E_0 \big (  |  a_{i,n} (t)|^{p/2}_{\mathbb B}| X_i|^{p/2}_{\mathbb B} \big )$. Since it is a convex function, for any $t\in [0,1]$,
\beq \label{convex}
G_{i,n} (t) \leq \max \big ( G_{i,n} (0) , G_{i,n} (1) \big ) \leq  b_{i,n}^{p/2} \, .
\eeq
Starting from \eqref{avantdernier}, using \eqref{convex} and the fact that
$ ( \E_0 (| X_n|^p_{\mathbb B} ) )^{2/p}\leq b_{n,n}$,
the inequality \eqref{keyinemoment} follows. $\lozenge$

\bigskip

\noindent {\bf Proof of Theorem \ref{inequality-moment-banach-martingale}.} The proof follows the lines of the proof of Proposition 2.1 in \cite{Rio99}. The only difference is that Inequality \eqref{majD2mart} is used to get his bound (2.1). For
the reader's convenience, let us give the main steps of the proof. For any $t \in [0,\infty[$, let $\varphi_n(t) = \Vert   | M_{n-1} +t d_n |_{{\mathbb B}} \Vert_p^p$.
Using Taylor's integral formula at order two together with Inequality \eqref{majD2mart},
we infer that
\[
\varphi_n(t) \leq \varphi_n(0) + c_p \Vert   | d_n |_{{\mathbb B}} \Vert_p^2 \int_0^t (t-s) \big ( \varphi(s))^{1-2/p} : = \phi_n(t) \, .
\]
Proceeding as at the top of page 150 in \cite{Rio99}, it follows that for any non-negative real $x$,
\[
\phi'_n(x) \leq \Vert   | d_n |_{{\mathbb B}} \Vert_p \sqrt{\frac{pc_p}{(p-1)}}  \sqrt{ \big ( \phi_n(x) \big )^{2-2/p} -  \big ( \varphi_n(0) \big )^{2-2/p}  } \, .
\]
Next, using  lemma 2.1 in \cite{Rio99} and the arguments following it, we derive
\[
\phi'_n(x) \leq \Vert   | d_n |_{{\mathbb B}} \Vert_p \sqrt{p c_p}\big ( \phi_n(x) \big )^{1-2/p} \sqrt{ \big ( \phi_n(x) \big )^{2/p} -  \big ( \varphi_n(0) \big )^{2/p}  } \, ,
\]
and then
\[
\big ( \phi_n(x) \big )^{2/p}\leq    \big ( \varphi_n(0) \big )^{2/p}  + p^{-1} c_p x^2  \Vert   | d_n |_{{\mathbb B}} \Vert^2_p   \, .
\]
Since $\varphi_n(x) \leq \phi_n(x)$, it follows that
\[
\Vert   | M_n |_{{\mathbb B}} \Vert_p^2 = \big (\varphi_n(1) \big)^{2/p} \leq \Vert   | M_{n-1} |_{{\mathbb B}} \Vert_p^2 +  p^{-1} c_p \Vert
|d_n|_{{\mathbb B}} \Vert_p^2 \, ,
\]
proving the theorem. $\lozenge$

\section{Hoeffding type inequalities for martingales}
\setcounter{equation}{0}
In the following corollary, we give an exponential inequality for
the deviation of the ${\mathbb L}^q$-norm of martingales.
\begin{Corollary} \label{corPinelis}
Let $q \in [2, \infty[$ and $({\mathcal X}, {\mathcal A}, \mu)$ a measure space. Let $(d_i)_{i \in {\mathbb N}}$ be a sequence of martingale differences with values in ${\mathbb L}^q ={\mathbb L}^q ({\mathcal X}, {\mathcal A}, \mu)$ (equipped with the norm $| \cdot |_q$) with respect to an increasing filtration $({\mathcal F}_i)_{i \in {\mathbb N}}$. Assume that for all $i \in {\mathbb N}$, there exists a positive real $b$ such that $\Vert |d_i |_q \Vert_{\infty} \leq b$. Let $M_n = \sum_{i=1}^n d_i$. For any positive integer $n$ and any positive real $x$, the following inequality holds
\beq \label{inePinelis}
{\mathbb P} \Big ( \max_{1 \leq k \leq n } | M_k |_q \geq x \Big ) \leq
 \begin{cases}
1 \quad \quad \quad \quad \quad \quad \quad \quad  \text{ if $x < b \sqrt{(q-1)n}$}\\
\frac{ (b^2 (q-1)n)^{q/2}}{x^q} \quad \quad \quad \ \,   \text{ if $ b \sqrt{(q-1)n}< x <   b \sqrt{ e(q-1) n}$} \\
\frac {1}{\sqrt e} \exp  \Big (  - \frac{x^2}{2 e b^2 n } \Big )
\quad  \text{ if $ x \geq   b \sqrt{ e(q-1) n}$} \, .
\end{cases}
\eeq
\end{Corollary}
Under the assumptions of Corollary \ref{corPinelis}, Theorem 3.5 in \cite{Pi} gives the following upper bound: for any positive integer $n $ and any positive real $x$,
\beq \label{inePinelis94}
{\mathbb P} \Big ( \max_{1 \leq k \leq n } | M_k |_q \geq x \Big ) \leq 2 \exp  \Big (  - \frac{x^2}{2 (q-1) b^2 n } \Big ) \, .
\eeq
It is noteworthy to indicate that for any $q \geq  e +1$, the bound in \eqref{inePinelis} is always better than the one given in \eqref{inePinelis94}.

\medskip

\noindent {\bf Proof of Corollary \ref{corPinelis}.} Let $p$ be a real number in $[2, \infty[$. By the Doob-Kolmogorov maximal inequality,
\[
{\mathbb P} \Big ( \max_{1 \leq k \leq n } | M_k |_q \geq x \Big ) \leq x^{-p} \E \big ( |M_n |_q^p \big ) \, .
\]
Therefore, using Inequality \eqref{IneThinemomentbismart}, we derive that for any $p \geq q$,
\[
{\mathbb P} \Big ( \max_{1 \leq k \leq n } | M_k |_q \geq x \Big ) \leq \left ( \frac{  \sqrt { a_p b^2 n}}{x}  \right )^p \, , \text{where $a_p =\max(p,q)-1$}\, .
\]
Taking $p=q$ if $x < \big (  (q-1) e b^2 n \big )^{1/2}$ (so in this case $a_p = q-1$) and $p= 1+  \frac{x^2}{ e b^2 n } $ if $x \geq \big (  (q-1) e b^2 n \big )^{1/2}$ (so in this case $a_p= p-1$), the inequality  \eqref{inePinelis} follows. $\lozenge$

\medskip

In the following corollary, we give an exponential inequality for
the deviation of the ${\mathbb L}^q$-norm of partial sums.
The proof  is omitted since it is exactly the same as that of Corollary \ref{corPinelis},
by using Inequality \eqref{IneThinemomentbis}
 instead of Inequality \eqref{IneThinemomentbismart}.

\begin{Corollary} \label{corPinelisgeneral}
Let $q \in [2, \infty[$ and $({\mathcal X}, {\mathcal A}, \mu)$ a measure space. Let $(X_i)_{i \in {\mathbb N}}$ be a sequence of random variables with values in ${\mathbb L}^q ={\mathbb L}^q ({\mathcal X}, {\mathcal A}, \mu)$ (equipped with the norm $| \cdot |_q$). Let $({\mathcal F}_i)_{i \geq 0}$ be an increasing sequence of $\sigma$-algebras such that $X_i$ is ${\mathcal F}_i$-measurable, and denote by
$\E_i ( \cdot) = \E ( \cdot | {\mathcal F}_i)$ the conditional expectation with
respect to ${\mathcal F}_i$.
For any positive integer $n$, let $S_n=\sum_{i=1}^n X_i$. Assume that for any integer $i \in [1,n]$,
$$
\Big \Vert |X_i|_q  \Big | \sum_{k=i}^n \E_i(X_k)  \Big |_q \Big  \Vert_{\infty} \leq b_n^2 \, .
$$
Then, for any positive real $x$, the following inequality holds
\begin{equation*} \label{inePinelisbis}
{\mathbb P} \Big ( | S_n |_q \geq x \Big ) \leq
\begin{cases}
1 \quad \quad \quad \quad \quad \quad \quad \quad \, \text{ if $x < b_n \sqrt{2(q-1)n}$}\\
\frac{ (2b_n^2 (q-1)n)^{q/2}}{x^q} \quad \quad \quad \ \,
\text{if $ b_n \sqrt{2(q-1)n}< x <   b_n \sqrt{2 e(q-1) n}$} \\
\frac{1}{\sqrt e}\exp  \Big (  - \frac{x^2}{4 e b_n^2 n } \Big ) \quad  \text{ if $ x \geq   b_n \sqrt{ 2e(q-1) n}$.}
\end{cases}
\end{equation*}
\end{Corollary}

\section{Moment and deviation inequalities for the empirical process of nonuniformly expanding maps}
\setcounter{equation}{0}
In this section, we shall apply Theorem \ref{inequality-moment-banach} and the inequalities recalled in Appendix
to obtain moment and deviation inequalities for the ${\mathbb L}^q$ norm of the centered
empirical distribution function of nonuniformly
expanding
 maps of the interval. More precisely, our results apply to the iterates of a map $T$ from $[0,1]$ to $[0,1]$ that can be modelled by a Young tower
with polynomial tails of the return time.

In Section \ref{YT}, we recall the formalism of Young towers, which has been described in many
papers (see for instance \cite{Young99} and \cite{MD}) with sometimes slight differences. Here we borrow the formalism described in Chapter 1 of Gou\"ezel's PhD thesis
\cite{Go02}.

The moment inequalities are stated  in Section \ref{modev}, and an application to the Wassertein metric between the empirical measure of $\{T, T^2, \ldots , T^n\}$ and the $T$-invariant distribution is presented in Section \ref{Wasser}. To be complete, we give in Section  \ref{simple}  some upper bounds for the maximum
of the partial sums of H\"older observables, which can be proved as in Section \ref{modev}.

\subsection{One dimensional maps modelled by Young towers}\label{YT}
Let $T$ be a map from $[0,1]$ to $[0,1]$, and $\lambda$ be a probability measure on $[0,1]$.
Let $Y$ be a Borel set of $[0,1]$, with $\lambda(Y)>0$.
Assume that there exist a partition (up to a negligible set)
$\{Y_k\}_{k \in \{1, \ldots, K\}}$ of $Y$ (note that $K$ can be infinite) and a sequence
$(\varphi_k)_{k \in \{1, \ldots, K\}}$ of increasing numbers such that
$T^{\varphi_k}(Y_k)=Y$. Let then $\varphi_Y$ be the function from
$Y$ to $\{\varphi_k\}_{k \in \{1, \ldots, K\}}$ such that
 $\varphi_Y(y)=\varphi_k$ if $y \in Y_k$.

 We then define a space
 $$
 X=\{(y,i): y \in Y, i < \varphi_Y(y)\}
 $$
 and a map $\bar T$ on $X$:
 $$
 \bar T(y,i)=\begin{cases}
  (y, i+1) \quad \quad \quad   \text{if $i< \varphi_Y(y)-1 $}\\
  (T^{\varphi_Y(y)}(y), 0) \quad  \text{if $i= \varphi_Y(y)-1$.}
  \end{cases}
 $$
 The space $X$ is the Young tower. One can define the floors $\Delta_{k,i}$ for
 $k \in \{1, \ldots, K\} $ and $i \in \{0, \ldots, \varphi_k -1\}$:
 $\Delta_{k,i}= \{ (y, i): y \in Y_k\}$.
 These floors define a partition of the tower:
 $$
 X= \bigcup_{k \in \{1, \ldots, K\}, i \in \{0, \dots , \varphi_k-1\}} \Delta_{k,i} \, .
 $$

 On $X$,  the measure $m$ is defined as follows:
 if $\bar B$ is a set included in $\Delta_{k,i}$, that can be written
 as $\bar B= B\times \{ i \}$ with $B \subset Y_k$,
 then $m(\bar B)=\lambda(B)$.
 Consequently, for a set $\bar A \subset \bigcup_{\{k \, : \, \varphi_k >i\} } \Delta_{k,i}$,
 which can be written as $\bar A= A \times \{ i \}=\big(\bigcup_{\{k \,  : \, \varphi_k >i\}} B_k \big) \times \{ i \}$
 with $B_k \subset Y_k$, one has
 $$
 m(\bar A)= \lambda(A)=
 \sum_{\{k \, : \, \varphi_k >i\} } \lambda (B_k).
 $$


 Let $\pi$ be the ``projection'' from $X$ to $[0,1]$ defined by
 $
 \pi(y,i)=T^i(y)
 $. Then, one has
 $$
 \pi \circ \bar T= T \circ \pi \, .
 $$
 Indeed, if $i < \varphi_Y(y)-1$, then
 $\bar T(y,i)=(y, i+1)$ so that
 $$
 \pi \circ \bar T(y,i)= \pi (y, i+1)= T^{i+1}(y)= T \circ \pi (y,i) \, .
 $$
 If $i= \varphi_Y(y)-1$, then
 $\bar T(y,i)=(T^{\varphi_Y(y)}(y), 0)$ so that
 $$
 \pi \circ \bar T(y,\varphi_Y(y)-1)=
  T^{\varphi_Y(y)}(y)= T(T^{\varphi_Y(y)-1}(y))=T \circ \pi (y,\varphi_Y(y)-1) \, .
 $$

 Assume now that  $\bar T$ preserves the probability $\bar \nu$ on $X$, and let
 $\nu$ be the image measure of $\bar \nu $ by $\pi$. Then, for any measurable and bounded function $f$,
 $$
 \nu(f(T))=\bar \nu (f(T \circ \pi))=\bar \nu ((f \circ \pi)(\bar T)) =
 \bar \nu ( f \circ \pi)= \nu(f) \, ,
 $$
 and consequently $\nu$ is invariant by $T$.

 The map $T$ can be modelled by a Young tower if:
 \begin{enumerate}
 \item For any $k \in \{1, \ldots, K\} $, $T^{\varphi_k}$ is a measurable isomorphism
 between $Y_k$ and $Y$. Moreover
 there exists  $C>0$ such that,
 for any   $k \in \{1, \ldots, K\} $ and almost every $x,y$ in $Y_k$,
 $$
   \Big |1-\frac{(T^{\varphi_k})'(x)}{(T^{\varphi_k})'(y)}\Big| \leq C
   |T^{\varphi_k}(x)-T^{\varphi_k}(y)| \, .
 $$
 \item There exists $C>0$ such that, for any   $k \in \{1, \ldots, K\} $ and
 almost every $x,y$ in $Y_k$, for any $i < \varphi_k$,
 $$
   |T^i(x)-T^i(y)| \leq C |T^{\varphi_k}(x)-T^{\varphi_k}(y)| \, .
 $$
 \item There exists $\tau>1$ such that, for any   $k \in \{1, \ldots, K\} $
 and almost every $x,y$ in $Y_k$:
 $$
   |T^{\varphi_k}(x)-T^{\varphi_k}(y)|\geq \tau |x-y| \, .
 $$
 \item $\sum_{k=1}^K \varphi_k \lambda(Y_k)< \infty$.
 \end{enumerate}
 If $T$ can be modelled by a Young tower, then,
 on the tower, there exists a unique $\bar T$-invariant probability measure $\bar \nu$
 which is absolutely continuous with respect to $m$.
 Hence, there exists a unique $T$-invariant
 measure
 $\nu$ which is absolutely continuous with respect to the measure $\lambda$ (see \cite{Go02}, Proposition 1.3.18). This
 measure is the image measure of $\bar \nu$ by the projection $\pi$ and is supported
 by
 $$\Lambda= \bigcup_{n\geq 0} T^n(Y). $$

 Let
$\bar Y$ be the basis of the tower, that is
$
  \bar Y=\{(y,0), y \in Y\}.
$
Let $\varphi_{\bar Y}$  be the function from
$\bar Y$ to $\{\varphi_k\}_{k \in \{1, \ldots, K\}}$ such that
$\varphi_{\bar Y}((y,0))=\varphi_{Y}(y)$.
By definition of $\bar T$ one gets
$\bar T^{\varphi_k}(\Delta_{k,0})=\bar Y$. In addition, the quantity $\bar \nu (\{(y, 0) \in \bar Y :
   \varphi_{\bar Y}((y,0))>k\})$ is exactly of the same order as $
   \lambda (\{y \in Y : \varphi_Y(y)>k\})$ (see \cite{Go02}, Proposition 1.1.24).

 On the tower, one defines the distance $s$ as follows:
 $s(x,y)=0$ is $x$ and $y$ do not belong to the same partition element
 $\Delta_{k,i}$. If $x=(a,i)$ and $y=(b,i)$ belong to the same $\Delta_{k,i}$
 (meaning that $a$ and $b$ belong to $Y_k$), then
 $\delta(x,y)= \beta^{s(x,y)}$ for  $\beta =1/\tau$, where
 $s(x,y)$ is the smallest integer $n$ such that $S^n(a)$ and $S^n(b)$
 are not in the same $Y_j$.

 Because of Item 3, we know that $|S'|\geq\tau>1$, so that
 $S$ is uniformly expanding. For $x=(a,i)$ and $y=(b,i)$ in $\Delta_{k,i}$, one has
 $$
   |\pi(x)-\pi(y)|= |T^i(a)-T^i(b)|\leq C |T^{\varphi_k}(a)-T^{\varphi_k}(b)|
 $$
 by Item 2. Since $T^{\varphi_k}=S$ on $Y_k$, and since $|S'|\geq \tau$, it follows that
 $$
  |\pi(x)-\pi(y)|\leq C\beta^{s(x,y)-1}\leq \frac C \beta \beta^{s(x,y)}\, .
 $$
 Now, if $x$ and $y$ do not belong to the same partition element
 $\Delta_{k,i}$, then $|\pi(x)-\pi(y)|\leq \beta^{s(x,y)}=1$. It follows
 that there exists a positive constant $K$ such that
 $$
 |\pi(x)-\pi(y)|\leq K \beta^{s(x,y)}\, ,
 $$
 meaning that $\pi$ is Lipschitz with respect to the distance $\delta$.

 Among the maps that can be modelled by a Young tower, we shall consider
 the maps defined as follows.
 \begin{Definition} \label{definitionYT}
 One says that the map $T$  can be modelled by a Young tower with polynomial tails of the return times of order $1/\gamma$ with $\gamma \in (0,1)$ if
$\lambda (\{ y \in Y :  \varphi_Y(y) >k\}) \leq Ck^{-1/\gamma}$.
\end{Definition}

Let us briefly describe some properties of such maps.
For $\alpha \in (0,1]$, let $\delta_\alpha= \delta^\alpha$, let
 $L_\alpha$ be the space of Lipschitz functions with respect to $\delta_\alpha$,
 and let
 \begin{equation}\label{seminorm}
 L_\alpha(f)= \sup_{x, y \in X}\frac{|f(x)-f(y)|}{\delta_\alpha(x,y)}.
 \end{equation}
 For any positive real $a$, let $L_{\alpha,a}$ be the set of functions such that $L_\alpha(f)\leq a$.

 Denote by $P$ the Perron-Frobenius operator of $\bar T$ with respect
 to $\bar \nu$: for any bounded measurable functions $\varphi, \psi,$
 $$
 \bar \nu (\varphi \cdot \psi \circ \bar T)= \bar \nu ( P(\varphi) \psi) \, .
 $$

Let $T$ be a map that can be modelled by a Young tower with polynomial tails of the return times of order $1/\gamma$.
Then one can prove that (see \cite{MD} and Lemma 2.2 in \cite{DP}):
for any $m \geq 1$ and any $\alpha \in (0,1]$, there exists $C_\alpha>0$ such that, for
 any $\psi \in L_\alpha$,
 \beq \label{contracPm}
  |P^m(\psi) (x) - P^m (\psi) (y) | \leq C_\alpha \delta_\alpha(x,y) L_\alpha(\psi)\, .
\eeq
Moreover, starting from the results by Gou\"ezel \cite{Go02}, we shall prove in Proposition \ref{taualpha}
of the appendix that, for any $ \alpha \in (0,1]$ there
 exists $K_\alpha>0$ such that
\beq  \label{majorationpourbeta}
   \bar \nu \Big (\sup_{f \in L_{\alpha,1}} |P^n(f) -\bar \nu (f)|
     \Big) \leq \frac{K_\alpha}{n^{(1-\gamma)/\gamma}} \, .
\eeq

A well known example of map which can be modelled by a Young tower with polynomial tails
of the return times is the  intermittent map $T_\gamma$ introduced   by Liverani {\it et al.} \cite{LSV}: for $\gamma \in (0,1)$,
\beq \label{LSVmap}
   T_\gamma(x)=
  \begin{cases}
  x(1+ 2^\gamma x^\gamma) \quad  \text{ if $x \in [0, 1/2[$}\\
  2x-1 \quad \quad \quad \ \  \text{if $x \in [1/2, 1]$;}
  \end{cases}
\eeq

 For this map, $\lambda$ is the Lebesgue measure on $[0,1]$ and one can take $Y=]1/2, 1]$. Let
 $x_0=1$, and define recursively $x_{n+1}= T_\gamma^{-1}(x_n) \cap [0, 1/2]$.
 One can prove that $x_n=\frac 1 2 (\gamma n )^{-1/\gamma}$. Let then
 $y_n=T_\gamma^{-1}(x_{n-1}) \cap ]1/2, 1]$. The $y_k$'s are built in such a way
 that $Y_k=]y_{k+1}, y_k]$ is the set of points $y$ in $Y$ for which
 $T_\gamma^k(Y_k)=Y$. One can verify, by controling explicitely the distortion,
 that the items 1,2 and 3 are satisfied with $\varphi_k=k$. Item 4 follows from
 the fact that
  $\sum_{k=1}^\infty k \lambda(Y_k)\leq C
 \sum_{k=1}^\infty k  k^{-(\gamma+1)/\gamma} < \infty$,
 since $\gamma \in (0,1)$.
 Moreover, one has
 $$
 \lambda (\{ y \in Y :  \varphi_Y(y) >k\})= \sum_{i= k+1}^\infty \lambda(Y_i) \leq
  Ck^{-1/\gamma} \, ,
 $$
 so that the tail of the return times is of order $1/\gamma$.

\subsection{Moment and deviation inequalities for the empirical process}\label{modev}

For any $q \in [2, \infty [$, let
\beq \label{defDnq}
D_{n,q}  = \Big (  \int_0^1 |G_n(t)|^q  dt \Big )^{1/q}\, ,
\eeq
where $G_n$ is defined by
\beq \label{defGn}
G_n (t) = \sum_{k=1}^n \big ( {\bf 1}_{T^k \leq t} - \nu ([0,t])\big ) \, , \ t \in [0,1] \, .
\eeq
Applying Lemma 1 in \cite{DM07}, we see that
\[
\frac 1 n D_{n,q} = \sup_{f \in W_{q',1}} \Big | \frac 1 n \sum_{k=1}^n
\big(f(T^k) -\nu(f)\big)\Big |\, ,
\]
where $q' = q/(q-1)$ and $W_{q',1}$ is the Sobolev ball
\beq \label{sobolev}
W_{q',1}= \left \{ f \ : \ \int_0^1 |f'(x)|^{q'}dx \leq 1  \right \} \, .
\eeq
Consequently, a moment inequality on $D_{n,q}$ provides a concentration inequality
of the empirical measure of $\{T, T^2, \cdots, T^n\}$ around $\nu$, on a class
of smooth functions. Note that, the class  $W_{q',1}$ is larger as $q$ increases, and always contains the
class of Lipschitz functions with Lipschitz constant 1.

In what follows, we shall denote by $\Vert \cdot \Vert_{p,\nu}$ the ${\mathbb L}^p$-norm on $([0,1],\nu)$
\begin{Theorem} \label{appYT1}
Let $T$ be a map that can be modelled by a Young tower with polynomial tails of the return times of order $1/\gamma$ with $\gamma \in (0,1/2)$, and let  $p_{\gamma}= 2 (1- \gamma)/ \gamma$.  For  $q \in [2, \infty [$ let $D_{n,q}$ be defined by \eqref{defDnq}. Then, there exists a positive constant $C$ such that for any $n \geq 1$,
\[
\Big \Vert \max_{1 \leq k \leq n} D_{k,q} \Big \Vert_{p_{ \gamma}, \nu}
\leq C \sqrt n \, .
\]
\end{Theorem}
As a consequence of Theorem \ref{appYT1}, for any $\gamma \in (0,1/2)$ and any positive real $x$,
\[
\nu \Big( \max_{1 \leq k \leq n} D_{k,q} \geq x \sqrt n
\Big) \leq  \frac{C}{x^{2 (1- \gamma)/ \gamma}} \, .
\]
In addition, proceeding as at the beginning of page 872 of the paper \cite{CG}, we infer that, under the assumptions of Theorem  \ref{appYT1}, for any real $p > 2 (1- \gamma)/ \gamma$, there exists a positive constant $C$ such that, for any $n \geq 1$,
\[
\Big \Vert \max_{1 \leq k \leq n} D_{k,q} \Big \Vert_{p,\nu} \leq C n^{(\gamma p + \gamma-1)/(\gamma p )} \, .
\]

Let us examine now the case where $\gamma \geq 1/2$.
\begin{Theorem} \label{appYT2}
Let $T$ be map that can be modelled by a  Young tower with polynomial tails of the return times of order $1/\gamma$ with $\gamma \in [1/2,1)$. For $q \in [2, \infty [$,  let $D_{n,q}$ be defined by \eqref{defDnq}.
\begin{enumerate}
\item There exists a positive constant $C$ such that for any $n \geq 1$,
$$
\Big \Vert \max_{1 \leq k \leq n} D_{n,q} \Big \Vert_{1/\gamma,\nu}
\leq C ( n \log n )^{\gamma}
\, .
$$
\item If $p > 1/\gamma$, then there exists a positive constant $C$ such that for any $n \geq 1$,
$$
\Big \Vert \max_{1 \leq k \leq n}D_{k,q} \Big \Vert_{p,\nu} \leq C n^{(\gamma p + \gamma-1)/(\gamma p )}\, .
$$
\end{enumerate}
\end{Theorem}

For the optimality of the moment bounds of Theorems  \ref{appYT1}
and \ref{appYT2}, we refer the paper by Melbourne and Nicol \cite{MN} and to the recent paper
by Gou\"ezel and Melbourne \cite{GoMe}. Since, for $q\geq 2$, the class $W_{q',1}$
contains the class of Lipschitz functions with Lipschitz constant 1, one
can apply Proposition 1.1 and 1.2 in \cite{GoMe}, showing that these bounds are optimal.
See also Remark \ref{remappYT3} below for more comments about the optimality.

\begin{Theorem} \label{appYT3}
Let $T$ be map that can be modelled by a  Young tower with polynomial tails of the return times of order $1/\gamma$ with $\gamma \in (1/2,1)$. For $q \in [2, \infty [$, let $D_{n,q}$ be defined by \eqref{defDnq}. Then, there exists a positive constant $C$ such that for any $n \geq 1$ and any positive real $x$,
\beq \label{resthmGPM3}
\nu \Big ( \max_{1 \leq k \leq n} D_{k,q} \geq x \, n^{\gamma} \Big ) \leq C x^{-1/\gamma}  \, .
\eeq
\end{Theorem}
Applying Theorem \ref{appYT3}, one gets for $p \in [1, 1/\gamma[$,
$$
\Big \Vert \max_{1 \leq k \leq n} D_{k,q} \Big \Vert_{p,\nu}^p
= p\int_0^\infty x^{p-1} \nu \Big ( \max_{1 \leq k \leq n} D_{k,q} \geq x  \Big ) dx
\leq p \int_0^{n^\gamma} x^{p-1} dx + C n p \int_{n^\gamma}^\infty  \frac{1}{x^{1+\gamma^{-1} -p}} dx \, .
$$
Consequently, for $p \in [1, 1/\gamma[$, there exists a positive constant $C$
such that
$$
\Big \Vert \max_{1 \leq k \leq n} D_{k,q} \Big \Vert_{p,\nu} \leq C n^\gamma \, .
$$

\begin{Remark} \label{remappYT3}  Inequality \eqref{resthmGPM3} cannot hold for $\gamma=1/2$. Indeed, for the map $T_{\gamma} $ defined in \eqref{LSVmap}, Item 1 of Theorem 1.1 in \cite{DDT} implies that, for any positive real $x$,
\[
\lim_{n \rightarrow \infty} \nu \Big (\frac{1}{\sqrt{ n \log n}} D_{n,2} >x  \Big ) = {\mathbb P} ( |N| >x) >0 \, ,  \]
where $N$ is a real-valued centered Gaussian random variable with positive variance. In addition, for $\gamma \in (1/2, 1)$, Item 2 of the same paper implies that
\[
\lim_{n \rightarrow \infty} \nu \Big (\frac{1}{n^{\gamma}} D_{n,2} >t  \Big ) =
{\mathbb P} ( |Z_{\gamma}| >t) >0 \, ,  \]
where $Z_{\gamma}$ is an $1/\gamma$-stable random variable such that $\lim_{x \rightarrow \infty} x^{1/\gamma} {\mathbb P} (|Z_{\gamma}| >x) =c >0$.
\end{Remark}

\subsection{Application to the Wasserstein metric between the empirical measure and the
invariant measure}\label{Wasser}
Let us give an application of the results of Section \ref{modev} to the Wasserstein distance between the empirical measure of $\{T, T^2,  \ldots, T^n\}$
and the invariant distribution $\nu$. Recall that Wasserstein distance $W_1$ between two probability measures $\nu_1$ and $\nu_2$ on $[0,1]$
is defined as
$$
  W_1(\nu_1,\nu_2)= \inf \Big \{ \int|x-y| \mu(dx,dy),
  \mu \in {\mathcal M}(\nu_1,\nu_2) \Big \}\, .
$$
where ${\mathcal M}(\nu_1,\nu_2)$ is the set of probability measures on
$[0,1]\times [0,1]$ with margins
$\nu_1$ and $\nu_2$.
Recall also that, in this one dimensional setting,
$$
W_1(\nu_1,\nu_2)= \int_0^1 |F_{\nu_1}(t)-F_{\nu_2}(t)| dt \, ,
$$
where
$F_{\nu_1}$ and $F_{\nu_2}$ are the distribution functions of $\nu_1$ and $\nu_2$ respectively. Therefore, setting
$$
  \nu_n= \frac{1}{n}\sum_{i=1}^n \delta_{T^i}
$$
we get that for any $q \geq 2$,
\[
W_1(\nu_n, \nu) \leq \frac 1 n  D_{n,q} \, .
\]
The following corollary is a direct consequence of the results of Section
\ref{modev}.
\begin{Corollary} \label{appKR}
Let $T$ be map that can be modelled by a  Young tower with polynomial tails of the return times of order $1/\gamma$ with $\gamma \in (0,1)$.
\begin{enumerate}
\item If $\gamma \in (0,1/2)$, then $ \Vert W_1(\nu_n, \nu) \Vert^p_{p,\nu} \ll n^{ - (1-\gamma)/\gamma  }$ for any $p \geq 2(1-\gamma)/\gamma$.
\item If $\gamma \in [1/2,1)$, then
\[
\Vert W_1(\nu_n, \nu) \Vert^p_{p,\nu}  \ll
\begin{cases}
n^{ - (1-\gamma)/\gamma  }  \log n \quad  \text{if $p=1/\gamma$}\\
n^{ - (1-\gamma)/\gamma  } \quad \quad \quad \,  \text{if $p>1/\gamma$.}
\end{cases}
\]
\item  If $\gamma \in (1/2,1)$, then, for any $n \geq 1$ and any positive real $x$,
\[
\nu \big ( W_1(\nu_n, \nu)  \geq x \, n^{\gamma-1} \big ) \ll x^{-1/\gamma}   \, .
\]
\end{enumerate}
\end{Corollary}

In their Theorem 1.4, Gou\"ezel and Melbourne \cite{GoMe} obtain  general bounds for
the moment of
separately Lipschitz functionals of $(T, T^2, \ldots, T^n)$, where $T$
is a (non necessarily one-dimensional) map  that can be modelled by a  Young tower with polynomial tails of the return times.

As a consequence of their results, one gets  the same inequalities as in
Corollary \ref{appKR}
but for the quantity $W_1(\nu_n, \nu)-{\mathbb E}(W_1(\mu_n, \nu))$ instead of
$W_1(\nu_n, \nu)$. Note that the upper bounds
 for $W_1(\nu_n, \nu)-{\mathbb E}(W_1(\mu_n, \nu))$  are valid  if $T$ is nonuniformly expanding from ${\mathcal X}$
to ${\mathcal X}$, where ${\mathcal X}$ can be any bounded metric space.

 The two results are not of the same nature. However, in our one dimensional setting,
the moments bounds of Corollary \ref{appKR} imply the same moment bounds
for $W_1(\nu_n, \nu)-{\mathbb E}(W_1(\mu_n, \nu))$,  because
$({\mathbb E}(W_1(\nu_n, \nu)))^p \leq \|W_1(\mu_n, \nu)\|_p^p$. The same remark does not hold for the deviation bounds, which are not directly comparable.

To conclude this section, let us mention that there is no hope to extend
Corollary \ref{appKR} to higher dimension with the same bounds. To see this, let us consider the  case of  ${\mathbb R}^d$-valued random variables $(X_1,X_2, \ldots, X_n)$
that are bounded,  independent, and identically distributed. Let $\nu_n$ be the empirical measure of
$\{X_1, X_2, \ldots, X_n \}$ and $\nu$ be the common distribution of the $X_i$'s.
It is well known that, when   $d\geq 3$ and $\nu$ has a component which is absolutely continuous with respect to the Lebesgue measure, ${\mathbb E}(W_1(\nu_n, \nu))$ is
exactly  of order $n^{-1/d}$, which is much
slower than $n^{-1/2}$.

\subsection{Moment and deviation inequalities for partial sums}\label{simple}
In this section, we  assume that $T$ is a nonuniformly expanding map on $({\mathcal X}, \lambda)$ with $\lambda$ a probability measure on ${\mathcal X}$, and  that $T$ can be modelled by a Young tower. Contrary to the previous sections, ${\mathcal X}$ can be any bounded  metric space and not necessarily the unit interval.
Let $f$ be a H\"older continuous function from ${\mathcal X}$ to ${\mathbb R}$ and $S_n(f) = \sum_{i=1}^{n} ( f \circ T^i - \nu (f) )$.

\begin{Theorem} \label{appMB}
Let $T$ be map that can be modelled by a  Young tower with polynomial tails of the return times of order $1/\gamma$ with $\gamma \in (0,1)$.
\begin{enumerate}
\item If $\gamma \in (0,1/2)$ then $\displaystyle \Big \Vert \max_{1\leq k \leq n} |S_k(f)| \Big \Vert^p_{p,\nu} \ll n^{p  - (1-\gamma)/\gamma  }$ for any $p \geq 2(1-\gamma)/\gamma$.
\item If $\gamma \in [1/2,1)$, then
\[
\Big \Vert \max_{1\leq k \leq n} |S_k(f)| \Big  \Vert^p_{p,\nu}  \ll
\begin{cases}
n \log n & \mbox{ if $p=1/\gamma$}\\
n^{p - (1-\gamma)/\gamma  } & \mbox{ if $p>1/\gamma$} \, .
\end{cases}
\]
\item  If $\gamma \in (1/2,1)$, for any $n \geq 1$ and any positive real $x$,
\[
\nu \Big ( \max_{1\leq k \leq n} |S_k(f)| \geq x \, n^{\gamma} \Big ) \ll x^{-1/\gamma}   \, .
\]
\end{enumerate}
\end{Theorem}

The proof is omitted since it is a simpler version of the proofs of Theorems \ref{appYT1}, \ref{appYT2} and \ref{appYT3}. Indeed the norm $| \cdot |_q$ is replaced by the absolute values and we do not need to deal with the supremum over a subset of the class of H\"older functions of order $1/q$.

After this paper was written, we became aware
that, using different methods based on martingale approximations, Gou\"ezel and Melbourne \cite{GoMe}
had independently obtained the upper bounds given in Theorem \ref{appMB}
(but for $|S_n(f)|$ instead of $\max_{1\leq k \leq n} |S_k(f)|$).

As in Section \ref{modev}, applying Propositions 1.1 and 1.2 in  \cite{GoMe},
we see that the moments bounds of Theorem \ref{appMB} cannot be improved.


Note also that, for the map $T_{\gamma} $ defined in \eqref{LSVmap},  we can make a similar remark as Remark \ref{remappYT3}:
Firstly,
Inequality \eqref{resthmGPM3} cannot hold for $\gamma=1/2$. Indeed by  Item 3 page 88 \cite{Go04}, if $f(0)\neq \nu(f)$, for any positive real $x$,
\[
\lim_{n \rightarrow \infty} \nu \Big (\frac{1}{\sqrt{ n \log n}} |S_n(f)| >x  \Big ) = {\mathbb P} ( |N| >x) >0 \, ,  \]
where $N$ is a real-valued centered Gaussian random variable with positive variance. In addition, for $\gamma \in (1/2, 1)$, Theorem 1.3 of the same paper implies that
\[
\lim_{n \rightarrow \infty} \nu \big ( |S_n(f)| > x  n^{\gamma} \big ) =
{\mathbb P} ( |Z_{\gamma}| >x) >0 \, ,  \]
where $Z_{\gamma}$ is an $1/\gamma$-stable random variable such that
$\lim_{x \rightarrow \infty} x^{1/\gamma} {\mathbb P} (|Z_{\gamma}| >x) =c >0$.

\subsection{Proofs of Theorems \ref{appYT1}, \ref{appYT2} and \ref{appYT3}.}
\noindent {\bf Proof of Theorem \ref{appYT1}.}  For any $t$, let $f_t $ be the function defined by $f_t ( x) = {\mathbf 1}_{x \leq t}$. Notice first that, for any $p \geq 1$,
\begin{multline*}
\Big \Vert \max_{1 \leq k \leq n} D_{k,q} \Big \Vert^p_{p,\nu} = \nu \Big ( \max_{1 \leq k \leq n} \Big | \int_0^1 \Big | \sum_{i=1}^k ({\mathbf 1}_{T^i \leq t} )  - \nu ([0,t])  \Big |^q dt \Big |^{p/q}\Big ) \\
= \bar{\nu} \Big ( \max_{1 \leq k \leq n} \Big | \int_0^1 \Big | \sum_{i=1}^k ( f_t \circ T^i \circ \pi  - \bar{\nu} (f_t \circ \pi) \Big  |^q dt \Big |^{p/q}\Big )\\ = \bar{\nu} \Big ( \max_{1 \leq k \leq n} \Big | \int_0^1 \Big  | \sum_{i=1}^k ( f_t  \circ \pi  \circ {\bar T}^i - \bar{\nu} (f_t \circ \pi)  \Big |^q dt \Big |^{p/q}\Big ) \, .
\end{multline*}
Let $g_t := f_t \circ \pi$ and $G(x) = \{ g_t(x) , t \in [0,1] \}$. Denote by $| \cdot |_q$  the norm associated to the Banach space ${\mathbb B}={\mathbb L}^q ( [0,1],dt)$. With these notations, we then have
\beq \label{p1th1MYT}
\Big \Vert \max_{1 \leq k \leq n} D_{k,q} \Big  \Vert^p_{p,\nu} = \bar{\nu} \Big ( \max_{1 \leq k \leq n} \Big |  \sum_{i=1}^k ( G({\bar T}^i) - \bar{\nu} (G({\bar T}^i)) )  \Big |_q^{p}\Big ) \, .
\eeq
Let now $(X_i)_{i \in {\mathbb N}}$ be a stationary Markov chain defined on a probability space $(\Omega, {\mathcal A}, {\mathbb P})$, with  state space $X$, transition probability $P$ and invariant distribution $\bar{\nu}$.  Recall then (see for instance Lemma XI.3 \cite{HH}) that for every $n\ge 1$, we have the following equalities  in law
(where in the left-hand side the law is meant under $\bar{\nu}$ and in the right-hand
side the law is meant under ${\mathbb P}$)
\begin{gather}
\label{law1bis} ({\bar T}^n,\ldots , {\bar T})\overset{d}=(X_1,\ldots , X_n) \nonumber \\
\label{law3bis} \max_{1 \leq k \leq n} \big |  \sum_{i=1}^k ( G({\bar T}^i) - \bar{\nu} (G({\bar T}^i)) )  \big |_q \overset{d}= \max_{1 \leq k \leq n}  \big | \sum_{i=k}^n (G (X_i) - \E (G(X_i)) )   \big |^p_q
\, .
\end{gather}
Therefore, starting from \eqref{p1th1MYT} and using \eqref{law3bis}, we infer that for any real $p \in [1, \infty[$,
\begin{multline} \label{p2th1MYT}
\Big \Vert \max_{1 \leq k \leq n} D_{k,q} \Big \Vert^p_{p,\nu} = \E \Big ( \max_{1 \leq k \leq n} \Big | \sum_{i=k}^n (G (X_i) - \E (G(X_i)) )   \Big |^p_q\Big )\\ \leq 2^p \E \Big ( \max_{1 \leq k \leq n} \Big | \sum_{i=1}^k (G (X_i) - \E (G(X_i)) )   \Big |^p_q\Big ) \, .
\end{multline}
Whence, Theorem \ref{appYT1} will follow if one can prove that there exists a positive constant $C$ such that for any $n \geq 1$,
\beq \label{aim1appYT1}
\E \Big ( \max_{1 \leq k \leq n} \Big | \sum_{i=1}^k (G (X_i) - \E (G(X_i)) )   \Big |^{\frac{2(1-\gamma)}{\gamma}}_q\Big ) \leq C n^{\frac{1-\gamma}{\gamma}} \, .
\eeq
With this aim, we shall apply the Rosenthal type inequality \eqref{inerosbanach} given in Appendix, with $p=2(1-\gamma)/\gamma$ (note that $p>2$ since $\gamma \in (0,1/2)$). Letting ${\mathcal F}_k = \sigma( X_i , i \leq k)$ and $G^{(0)} = G - \E (G(X_1))$, this leads to
\begin{multline} \label{p3th1MYT}
\E \Big ( \max_{1 \leq k \leq n} \Big | \sum_{i=1}^k (G (X_i) - \E (G(X_i)) )   \Big |^{\frac{2(1-\gamma)}{\gamma}}_q\Big )
\ll n{\mathbb{E}}
\Big(|G (X_1)|_{q}^{\frac{2(1-\gamma)}{\gamma}}\Big)\\
+n\Big (   \sum_{k=1}^{n}\frac{1}{k^{1+\delta \gamma/(1-\gamma)}}\Big \Vert\E_0 \Big ( \Big  | \sum_{i=1}^k G^{(0)}(X_i)   \Big |^2_q \Big ) \Big  \Vert_{(1-\gamma)/\gamma}^{\delta}\Big )^{\frac{(1-\gamma)}{\delta \gamma}}\, .
\end{multline}
where $\delta = \min ( 1/2 , \gamma / (2-4\gamma) ) $. To handle the terms $\big \Vert\E_0 \big ( \big  | \sum_{i=1}^k G^{(0)}(X_i)   \big |^2_q \big ) \big \Vert_{(1-\gamma)/\gamma}$ in Inequality \eqref{p3th1MYT}, we shall use Inequality \eqref{IneThinemomentbis} which together with Item 1 of Lemma \ref{constantecp} leads to
\begin{multline*}
\E_0 \Big ( \Big  | \sum_{i=1}^k G^{(0)}(X_i)   \Big |^2_q \Big ) \leq 2 (2q-3) \sum_{i=1}^k \sum_{\ell=i}^k \E_0 ( |G^{(0)}(X_i)|_q | \E_i ( G^{(0)}(X_{\ell}) )|_q ) \\
\leq 2 (2q-3) \sum_{i=1}^k \sum_{\ell=i}^k \E_0 (  |  \E_i (  G^{(0)}(X_{\ell}) )|_q ) \, ,
\end{multline*}
where for the last inequality, we have used the fact that for any $i$, $ |G^{(0)}(X_i)|_q \leq 1$ almost surely. Hence
\beq\label{p3th1MYTbis}
\Big \Vert \E_0 \Big ( \Big  | \sum_{i=1}^k G^{(0)}(X_i)   \Big |^2_q \Big )\Big \Vert_{(1-\gamma)/\gamma}\leq 2 (2q-3) \sum_{i=1}^k \sum_{\ell=i}^k \Vert \E_0 (  |  \E_i (  G^{(0)}(X_{\ell}) )|_q )\Vert_{(1-\gamma)/\gamma} \, .
\eeq
Let us now handle the term $\Vert \E_0 (  |  \E_i (  G^{(0)}(X_{\ell}) )|_q )\Vert_{(1-\gamma)/\gamma}$ in Inequality \eqref{p3th1MYTbis}. With this aim, we first notice that
\[
|  \E_i ( G^{(0)}(X_{\ell}) )|^q_q  = \int_0^1 \big | \E ( {\mathbf 1}_{\pi (X_\ell) \leq t }  | X_i) - \E ( {\mathbf 1}_{\pi (X_\ell) \leq t } ) \big |^q dt
\]
Using Lemma 1 in \cite{DM07}, we have
\[
\int_0^1 \big | \E ( {\mathbf 1}_{\pi (X_\ell) \leq t }  | X_i) - \E ( {\mathbf 1}_{\pi (X_\ell) \leq t } ) \big |^q dt = \sup_{h \in W_{q',1}} \big | P_{\pi (X_{\ell})  | X_i} (h) - P_{\pi (X_{\ell})}(h) \big |^q  \, ,
\]
where  the Sobolev ball $W_{q',1}$ is defined in \eqref{sobolev},
 $P_{\pi (X_{\ell})  | X_i}$ is the conditional distribution
of $\pi (X_{\ell})$ given   $X_i$, and $P_{\pi (X_{\ell})}$ is the distribution of
$\pi(X_{\ell})$.
 Therefore
\begin{multline*}
|  \E_i ( G^{(0)}(X_{\ell}) )|_q  = \sup_{h \in W_{q',1}} \big | P_{\pi (X_{\ell})  | X_i} (h) - P_{\pi (X_{\ell})}(h) \big |
= \sup_{h \in W_{q',1}} \big | P_{X_{\ell}  | X_i} (h \circ \pi) - P_{X_{\ell}}(h \circ \pi) \big | \, ,
\end{multline*}
where $P_{X_{\ell}  | X_i}$ is the conditional distribution
of $X_{\ell}$ given   $X_i$, and $P_{X_{\ell}}$ is the distribution of
$X_{\ell}$.
Notice now that if $f \in  W_{q',1}$ then for any $x$ and $y$ in $[0,1]$,
\[ |f(x) - f(y) | = \Big | \int_x^y f'(t) dt  \Big | \leq | x-y|^{1/q}\Big (  \int_0^1 | f'(x)|^{q'}dx\Big )^{1/q'} \, . \]Therefore,
$$
W_{q',1} \subset H_{1/q,1} \, ,
$$
where $H_{1/q, 1}$ is the set of functions that are $1/q$-H\"older with
H\"older constant 1. It follows that, for any $h \in W_{q',1}$, there exists a positive constant $C$ such that
$$
| h \circ \pi (x) -  h \circ \pi (y) | \leq | \pi (x) - \pi (y) |^{1/q} \leq C \delta_{1/q} (x,y) \, ,
$$
proving that $ h\circ \pi$ belongs to the set $L_{1/q,C}$ defined right after
\eqref{seminorm}. Let now
\[
f_{\ell-i, h} (x) : =\big | P_{X_{\ell}  | X_i=x} (h \circ \pi) -  P_{X_{\ell}}(h \circ \pi)  \big | = \big |  P^{\ell-i} ( h \circ \pi) (x) - \bar{\nu} ( h \circ \pi) \big | \, .
 \]
Using the triangle inequality, we have
 \[ | f_{\ell-i, h} (x) - f_{\ell-i, h} (y) | \leq  \big |  P^{\ell-i} ( h \circ \pi) (x) - P^{\ell-i} ( h \circ \pi) (y)\big | \, .\]
Since $h \circ \pi$ belongs to $L_{1/q,C}$, the contraction property \eqref{contracPm} entails that
 \[
  | f_{\ell-i, h} (x) - f_{\ell-i, h} (y) | \leq C  C_{1/q} \delta_{1/q}(x,y)  \, .\]
Let $\widetilde C= C  C_{1/q}$.
We have shown that, for any $h \in W_{q',1}$, $f_{\ell-i,h} \in {\mathcal F}_{\ell-i} \subset L_{1/q,{\widetilde C}}$.
Then, setting $$m_{\ell-i}(x) = \sup_{h \in W_{q',1}} f_{\ell-i, h} (x) \, $$ we have $ m_{\ell-i}(x) = \sup_{g \in {\mathcal F}_{\ell-i}} g (x)$. Therefore, if $m_{\ell-i}(x) \geq m_{\ell-i}(y)$,
\[
m_{\ell-i}(x) - m_{\ell-i}(y) = g_x(x) - g_y(y) \leq g_x(x) - g_x(y) \leq {\widetilde C} \delta_{1/q} (x,y) \, ,
\]
since $ {\mathcal F}_{\ell-i} \subset L_{1/q,{\widetilde C}}$.
So overall,
\[
|  \E_i ( G^{(0)}(X_{\ell}) )|_q - \E |  \E_i ( G^{(0)}(X_{\ell}) )|_q = m_{\ell-i}(X_i)  - \E ( m_{\ell-i} (X_i) ) \, ,
\]
with $m_{\ell-i} \in L_{1/q, \widetilde C}$. Next, using \eqref{majorationpourbeta}, it follows that there exists a positive constant $C$ such that, for any $i \geq 1$,
\beq \label{p4th1MYT}
 \Vert  \E_0 \big( | \E_i ( G^{(0)}(X_{\ell}) )|_q \big )- \E |  \E_i ( G^{(0)}(X_{\ell}) )|_q  \Vert_1 =  \Vert P^i ( m_{\ell-i} )   - \bar{\nu}( m_{\ell-i} ) \Vert_1 \leq C i^{-(1-\gamma)/\gamma} \, .
\eeq
Using similar arguments we infer that there exists a positive constant $C$ such that, for any $\ell \geq i+1$,
\beq \label{p5th1MYT}
\Vert |  \E_i ( G^{(0)}(X_{\ell}) )|_q  \Vert_1
 = \Vert |  \E_0 ( G^{(0)}(X_{\ell-i}) )|_q   \Vert_1
 \leq \bar{\nu} \Big (  \sup_{g \in L_{1/q,{\widetilde C}}} \big | P^{\ell -i} (g) - \bar{\nu} (g) \big | \Big )  \leq C (\ell - i)^{-(1-\gamma)/\gamma} \, .
\eeq
We control now the quantity $\sum_{i=1}^k \sum_{\ell=i}^k \Vert \E_0 (  |  \E_i (  G^{(0)}(X_{\ell}) )|_q ) \Vert_{(1-\gamma)/\gamma}$ with the help of
\eqref{p4th1MYT} and \eqref{p5th1MYT}. With this aim, we first write the following decomposition:
\begin{multline*}
\sum_{i=1}^k \sum_{\ell=i}^k \Vert \E_0 (  |  \E_i (  G^{(0)}(X_{\ell}) )|_q ) \Vert_{(1-\gamma)/\gamma} \leq \sum_{i=1}^k \sum_{\ell=2i+1}^k \Vert  |  \E_i (  G^{(0)}(X_{\ell}) )|_q  \Vert_{(1-\gamma)/\gamma}\\ + \sum_{i=1}^k \sum_{\ell=i}^{2i} \Vert \E_0 (  |  \E_i (  G^{(0)}(X_{\ell}) )|_q ) - \E |  \E_i ( G^{(0)}(X_{\ell}) )|_q \Vert_{(1-\gamma)/\gamma} + \sum_{i=1}^k \sum_{\ell=i}^{2i} \Vert |  \E_i (  G^{(0)}(X_{\ell}) )|_q \Vert_1
\end{multline*}
Next, since $(1-\gamma)/\gamma >1$ and for any $i$, $|G^{(0)}(X_{i}) )|_q \leq 1$
almost surely,  we get
\begin{multline*}
\sum_{i=1}^k \sum_{\ell=i}^k \Vert \E_0 (  |  \E_i (  G^{(0)}(X_{\ell}) )|_q ) \Vert_{(1-\gamma)/\gamma} \leq \sum_{i=1}^k \sum_{\ell=2i+1}^k \Vert  |  \E_i (  G^{(0)}(X_{\ell}) )|_q  \Vert_1^{\gamma/(1-\gamma)}\\ + 2^{\frac{1-2\gamma}{1-\gamma}}\sum_{i=1}^k \sum_{\ell=i}^{2i} \Vert \E_0 (  |  \E_i (  G^{(0)}(X_{\ell}) )|_q ) - \E |  \E_i ( G^{(0)}(X_{\ell}) )|_q \Vert_1^{\gamma/(1-\gamma)} + \sum_{i=1}^k \sum_{\ell=i}^{2i} \Vert |  \E_i (  G^{(0)}(X_{\ell}) )|_q \Vert_1 \, .
\end{multline*}
Therefore, using \eqref{p4th1MYT} and \eqref{p5th1MYT}, we derive that
\begin{multline} \label{p6th1MYT}
\sum_{i=1}^k \sum_{\ell=i}^k \Vert \E_0 (  |  \E_i (  G^{(0)}(X_{\ell}) )|_q ) \Vert_{(1-\gamma)/\gamma}\\
\ll \sum_{i=1}^k \sum_{\ell=2i+1}^k \frac{1}{\ell -i} + \sum_{i=1}^k \sum_{\ell=i}^{2i} \frac{1}{i} + k+ \sum_{i=1}^k \sum_{\ell=i+1}^{2i} \frac{1}{(\ell - i)^{\frac{1-\gamma}{\gamma}}} \ll k\, .
\end{multline}
So starting from  \eqref{p3th1MYT} and taking into account \eqref{p3th1MYTbis},  \eqref{p6th1MYT} and the fact that $\gamma/(1-\gamma)< 1$, we get
\[
\E \Big ( \max_{1 \leq k \leq n} \Big | \sum_{i=1}^k (G (X_i) - \E (G(X_i)) )   \Big |^{\frac{2(1-\gamma)}{\gamma}}_q\Big )
\ll n+n\Big (   \sum_{k=1}^{n}\frac{k^{\delta}}{k^{1+\delta \gamma/(1-\gamma)}} \Big )^{(1-\gamma)/(\delta \gamma)} \ll n^{(1-\gamma)/\gamma}\, ,
\]
which completes the proof of \eqref{aim1appYT1} and then of the theorem. $\lozenge$

\bigskip

\noindent {\bf Proof of Theorem \ref{appYT2}.} We keep the same notations as in the proof of Theorem \ref{appYT1}.

We start by proving Item 1. By \eqref{p2th1MYT}, it suffices to prove that
there exists a positive constant $C$ such that for any $n \geq 1$,
\beq \label{aim1appYT2}
\E \Big ( \max_{1 \leq k \leq n} \Big | \sum_{i=1}^k (G (X_i) - \E (G(X_i)) )   \Big |^{1/ \gamma}_q \Big ) \leq C  \, n \log n\, .
\eeq
Assume first that $\gamma=1/2$. Applying Inequality \eqref{IneThinemomentbismax}, taking into account the stationarity and the fact that $ | G (X_1) - \E (G(X_1)) |_q \leq 1$ almost surely, we derive
\[
\E \Big ( \max_{1 \leq k \leq n} \big | \sum_{i=1}^k (G (X_i) - \E (G(X_i)) )   \big |^{1/ \gamma}_q \big )  \ll n + n\sum_{k=1}^n\Vert |  \E_0 ( G^{(0)}(X_{k}) )|_q  \Vert_1\, .
\]
Therefore, using \eqref{p5th1MYT}, it follows that
\[
\E \Big ( \max_{1 \leq k \leq n} \Big | \sum_{i=1}^k (G (X_i) - \E (G(X_i)) )   \Big |^{1/ \gamma}_q \Big )  \ll n + n\sum_{k=1}^n k^{-1}\, .
\]
proving \eqref{aim1appYT2} in the case $\gamma=1/2$. We turn now to the proof of \eqref{aim1appYT2} when $\gamma \in (1/2,1)$. With this aim, we apply the moment inequality (with $p = 1/\gamma$) stated in Proposition \ref{prop:maximal}. This leads to
\[
\E \Big ( \max_{1 \leq k \leq n} \Big | \sum_{i=1}^k (G (X_i) - \E (G(X_i)) )   \Big |^{1/ \gamma}_q \Big )  \leq C_{\gamma}n \sum_{k=0}^{n-1} (k+1)^{(1-2\gamma)/\gamma} \Vert |  \E_0 ( G^{(0)}(X_{k}) )|_q  \Vert_1 \, ,
\]
where $C_{\gamma}$ is a positive constant depending only on $ \gamma$. Therefore, for any $ \gamma \in (1/2,1)$ using \eqref{p5th1MYT}, we get
\[
\E \Big ( \max_{1 \leq k \leq n} \Big | \sum_{i=1}^k (G (X_i) - \E (G(X_i)) )   \Big |^{1/ \gamma}_q \Big )  \leq {\tilde C}_{\gamma} n \Big ( 1+ \sum_{k=1}^{n-1} k^{-1} \Big ) \, ,
\]
proving \eqref{aim1appYT2} in case $\gamma \in (1/2,1)$. This ends the proof of Item 1.

\smallskip

We turn now to the proof of Item 2. By  \eqref{p2th1MYT}, it suffices to prove that, for $\gamma \in [1/2,1)$ and  $p> 1/\gamma$,
there exists a positive constant $C$ such that for any $n \geq 1$,
\beq \label{aim1Th2YTitem2}
\E \Big ( \max_{1 \leq k \leq n} \Big | \sum_{i=1}^k (G (X_i) - \E (G(X_i)) )   \Big |^{p}_q \Big ) \leq C \, n^{p + (\gamma-1)/\gamma}\, .
\eeq We shall distinguish two cases: ($p\geq 2$ and $p>1/\gamma$) or $p\in ]1/\gamma , 2[$. We first consider the case where $p\geq 2$ and $p>1/\gamma$. To prove \eqref{aim1Th2YTitem2}, we shall apply Inequality \eqref{IneThinemomentbismax}.  Taking into account the stationarity and the fact that $ | G (X_1) - \E (G(X_1)) |_q \leq 1$ almost surely, we derive
\[
\E \Big ( \max_{1 \leq k \leq n} \Big | \sum_{i=1}^k (G (X_i) - \E (G(X_i)) )   \Big |^p_q \Big )   \ll n^{p/2} \Big (  \sum_{k=0}^n  \Vert |  \E_0 ( G^{(0)}(X_{k}) )|_q  \Vert_1^{2/p}\Big )^{p/2}  \, .
\]
Next, using \eqref{p5th1MYT} and the fact that $2(1-\gamma)/(\gamma p) <1$, Inequality \eqref{aim1Th2YTitem2} follows.

We consider now the case where $p\in ]1/\gamma , 2[$. Using, once again, the moment inequality stated in Proposition \ref{prop:maximal}, we get
\[
\E \Big ( \max_{1 \leq k \leq n} \Big | \sum_{i=1}^k (G (X_i) - \E (G(X_i)) )   \Big |^p_q \Big )  \leq C_{p}n \sum_{k=0}^{n-1} (k+1)^{p-2} \Vert |  \E_0 ( G^{(0)}(X_{k}) )|_q  \Vert_1\, ,
\]
where $C_{p}$ is a positive constant depending only on $ p$. Using then \eqref{p5th1MYT}  and the fact that $p>1/\gamma$, \eqref{aim1Th2YTitem2} follows. This ends the proof of the theorem. $\lozenge$

\bigskip

\noindent {\bf Proof of Theorem \ref{appYT3}.} We keep the same notations as in the proof of Theorem \ref{appYT1}. Notice first that, for any non-negative $x$,
\begin{align*}
\nu \Big (  \max_{1 \leq k \leq n} D_{k,q} \geq x \Big )
&= \bar{\nu} \Big ( \max_{1 \leq k \leq n} \Big | \int_0^1 \Big | \sum_{i=1}^k ( f_t \circ T^i \circ \pi  - \bar{\nu} (f_t \circ \pi) \Big  |^q dt \Big |^{1/q} \geq x \Big ) \\ &= \bar{\nu} \Big ( \max_{1 \leq k \leq n} \Big | \int_0^1 \Big  | \sum_{i=1}^k ( f_t  \circ \pi  \circ {\bar T}^i - \bar{\nu} (f_t \circ \pi)  \Big |^q dt \Big |^{1/q} \geq x \Big ) \\ &= \bar{\nu} \Big ( \max_{1 \leq k \leq n} \Big |  \sum_{i=1}^k ( G({\bar T}^i) - \bar{\nu} (G({\bar T}^i)) )  \Big |_q \geq x \Big ) \, .
\end{align*}
According to \eqref{law3bis},
\begin{multline*}
\nu \Big (  \max_{1 \leq k \leq n} D_{k,q} \geq x \Big )
= {\mathbb P} \Big (   \max_{1 \leq k \leq n} \Big | \sum_{i=k}^n (G (X_i) - \E (G(X_i)) )   \Big |_q \geq  x \Big ) \\
\leq {\mathbb P} \Big (   \max_{1 \leq k \leq n} \Big | \sum_{i=1}^k (G (X_i) - \E (G(X_i)) )   \Big |_q \geq  x/2 \Big ) \, .
\end{multline*}
The theorem will then follow if we can prove that, for any positive real $x$,
\beq \label{aimth3YT}
{\mathbb P} \Big (   \max_{1 \leq k \leq n} \Big | \sum_{i=1}^k (G (X_i) - \E (G(X_i)) )   \Big |_q \geq 4 x \Big ) \ll n x^{-1/\gamma} \, .
\eeq
To prove this inequality, we shall apply Proposition \ref{prop:maximal} with lag $[x]$. Using \eqref{p5th1MYT}, this leads to the following inequality: for any positive real $x$,
\[
{\mathbb P} \Big (   \max_{1 \leq k \leq n} \Big | \sum_{i=1}^k (G (X_i) - \E (G(X_i)) )   \Big |_q \geq 4 x \Big ) \ll \frac{n}{x^{1/\gamma}}+ \frac{n}{x^2}  \sum_{k=0}^{[x]}  \frac{1}{(k+1)^{(1-\gamma)/\gamma}}\, ,
\]
and \eqref{aimth3YT} follows. $\lozenge$

\section{Appendix}
\setcounter{equation}{0}

\subsection{A Rosenthal-type  inequality for stationary sequences}

In this section, for the reader convenience, we recall the Rosenthal-type inequality stated in \cite{DMP} (see Inequality (3.11) therein). This inequality is the extension to Banach-valued random variables of the Rosenthal type inequality given by Merlev\`ede and Peligrad \cite{MP}.

Let $(\Omega,{\mathcal A}, {\mathbb P} )$ be a
probability space, and $\theta:\Omega \mapsto \Omega$ be
a bijective bimeasurable transformation preserving the probability ${\mathbb P}$.
For a $\sigma$-algebra ${\mathcal F}_0$ satisfying ${\mathcal F}_0
\subseteq T^{-1 }({\mathcal F}_0)$, we define the nondecreasing
filtration $({\mathcal F}_i)_{i \in {\mathbb Z}}$ by ${\mathcal F}_i =\theta^{-i
}({\mathcal F}_0)$.
We shall use the notations $\E_k (\cdot) =
\E (\cdot | {\mathcal F}_k)$.

Let $X_0$ be a random variable with values in ${\mathbb B}$.
Define the stationary sequence  $(X_i)_{i \in \mathbb Z}$ by
$X_i = X_0 \circ T^i$, and the partial sum $S_n$ by $S_n=X_1+X_2+\cdots + X_n$.

\begin{Theorem} \label{Rosenthalbanachvalued}
Assume that $X_0$ belongs to ${\mathbb L}^p ({\mathbb B} ) $ where $({\mathbb B}, |\cdot|_{\mathbb B})$ is a separable Banach space and $p$ is a real number in $]2, \infty[$. Assume that $X_0$ is ${\mathcal F}_0 $-measurable. Then, for any $r \geq 0$,
\begin{equation} \label{rosenthaline1}
{\mathbb{E}} \Big (\max_{1\leq j\leq 2^r}|S_{j}|_{\mathbb B}^{p}\Big )\ll 2^r{\mathbb{E}
}(|X_{0}|^p_{\mathbb B})+ 2^r\left(  \sum_{k=0}^{r-1} \frac{
     \Vert{\mathbb{E}}_{0}( | S_{2^k} |_{\mathbb B}^{2})\Vert_{p/2}^{\delta}}{2^{2\delta k/ p}}\right)^{p/(2\delta)}\,,
\end{equation}
where $\delta=\min(1/2,1/(p-2))$.
\end{Theorem}
\begin{Remark} \label{remarknodiadic}
The inequality in the above theorem implies that for any positive integer $n$,
\beq \label{inerosbanach}
{\mathbb{E}} \Big (\max_{1\leq j\leq n}|S_{j}|_{\mathbb B}^{p}\Big )\ll n{\mathbb{E}
}(|X_{0}|_{\mathbb B})^{p}+n\left(  \sum_{k=1}^{n}\frac{1}{k^{1+2\delta
/p}} \Vert{\mathbb{E}}_{0}( | S_{k} |_{\mathbb B}^{2})\Vert_{p/2}^{\delta}\right)
^{p/(2\delta)}\, .
\eeq
\end{Remark}

\subsection{A deviation inequality}

The following proposition is adapted from Proposition 4 in  \cite{DMtpa}. It also extends Proposition 6.1 in \cite{DDT} to random variables taking values in a separable Banach space  belonging to the class ${\widetilde{\mathcal C}}_2 (2, {\tilde c}_2)$.

\begin{Proposition}\label{prop:maximal}
Let $Y_1, Y_2, \ldots, Y_n$ be $n$ random variables with values in a separable Banach
space $({\mathbb B}, |\cdot |_{\mathbb B})$ belonging to the class ${\widetilde{\mathcal C}}_2 (2, {\tilde c}_2)$. Assume that
 ${\mathbb P}(|Y_k|_{\mathbb B}\leq M)=1$ for
any $k \in \{1, \ldots , n\}$. Let ${\mathcal F}_1, \ldots ,
{\mathcal F}_n$ be an increasing filtration such that $Y_k$ is ${\mathcal F}_k$-measurable for any $k \in \{1, \ldots , n\}$.
Let $S_n=\sum_{k=1}^n Y_k$, and for
$k \in \{0, \ldots, n-1 \}$, let
\beq \label{defthetakappendix}
\theta(k)= \max \Big \{
 {\mathbb E}(|{\mathbb E}(Y_{i}|{\mathcal F}_{i-k})|_{\mathbb B}), i  \in \{k+1, \ldots , n\} \Big \} \, .
\eeq
Then, for any $q \in \{1, \ldots, n\}$, and any $x\geq qM$, the following inequality holds
\beq \label{propmaxIne1}
{\mathbb P} \Big ( \max_{1 \leq k \leq n} |S_k|_{\mathbb B} \geq 4x\Big)
\leq  \frac{n \theta(q)}{x}{\bf 1}_{q<n} + \frac{4{\tilde c}_2 K^2 nM}{x^2} \sum_{k=0}^{q-1} \theta (k)  \, ,
\eeq
where $K=\sqrt{\max ({\tilde c}_2 ,1)}$. In addition, for any $p \in [1,2[$,
\beq \label{propmaxIne2}
{\mathbb E} \big ( \max_{1 \leq k \leq n} |S_k|^p_{\mathbb B} \big ) \leq \Big ( 4^p p + \frac{4^{p+1} p {\tilde c}_2 K^2}{2-p} \Big ) M^{p-1} n \sum_{k=0}^{n-1} (k+1)\theta(k) \, .
\eeq
\end{Proposition}
\noindent {\bf Proof of Proposition \ref{prop:maximal}.} Let $S_0=0$ and define the random variables $U_i$ by:
 $U_i=S_{iq}-S_{(i-1)q}$ for
$i \in \{ 1, \ldots, [n/q]\}$ and $U_{[n/q]+1}=S_n-S_{q[n/q]}$. By Proposition 4 in  \cite{DMtpa}, for any $x\geq Mq$,
\begin{align}\label{propFlo}
{\mathbb P} \Big ( \max_{1 \leq k \leq n} |S_k|_{\mathbb B} \geq 4x\Big)
&\leq \frac{1}{x} \sum_{i=3}^{[n/q]+1}
{\mathbb E}( |{\mathbb E}(U_i|{\mathcal F}_{(i-2)q})|_{\mathbb B}) +
\frac{{\tilde c}_2}{x^2} \sum_{i=1}^{[n/q]+1}{\mathbb E}(|U_i-{\mathbb E}(U_i|{\mathcal F}_{(i-2)q}) |_{\mathbb B}^2) \nonumber  \\
&\leq
 \frac{1}{x} \sum_{i=3}^{[n/q]+1}
{\mathbb E}( |{\mathbb E}(U_i|{\mathcal F}_{(i-2)q})|_{\mathbb B}) +
\frac{4 {\tilde c}_2}{x^2} \sum_{i=1}^{[n/q]+1}{\mathbb E}(|U_i|_{\mathbb B}^2) \, .
\end{align}
Since $(\theta(k))_{k \geq 0}$ is a non-increasing sequence, it is not hard to see that
\begin{equation}\label{B2bis}
\sum_{i=3}^{[n/q]+1}
{\mathbb E}( |{\mathbb E}(U_i|{\mathcal F}_{(i-2)q})|_{\mathbb B})
\leq n \theta(q) {\bf 1}_{q<n}  \, .
\end{equation}
To handle the second term in (\ref{propFlo}), we use Inequality \eqref{IneThinemomentbis} with $p=2$. This leads to the following upper bounds: for any $i \in \{ 1, \ldots, [n/q]\}$,
\[
{\mathbb E}(|U_i|_{\mathbb B}^2)  \leq K^2 \sum_{k=(i-1)q+1}^{iq} \sum_{j=k}^{iq} \E \big ( | Y_k|_{\mathbb B} | \E ( Y_j | {\mathcal F}_k ) |_{\mathbb B} \big ) \, ,
\]
and
\[
{\mathbb E}(|U_{[n/q]+1}|_{\mathbb B}^2)  \leq K^2 \sum_{k=q[n/q] +1}^{n} \sum_{j=k}^{n} \E \big ( | Y_k|_{\mathbb B} | \E ( Y_j | {\mathcal F}_k ) |_{\mathbb B} \big ) \, ,
\]
where $K=\sqrt{\max ({\tilde c}_2 ,1)}$. Using the fact that ${\mathbb P}(|Y_k|_{\mathbb B}\leq M)=1$ for any $k \in \{1, \ldots , n\}$ and that $(\theta(k))_{k \geq 0}$ is a non-increasing sequence, we then derive that, for any $i \in \{ 1, \ldots, [n/q]\}$,
\[
{\mathbb E}(|U_i|_{\mathbb B}^2)  \leq K^2 M \sum_{k=(i-1)q+1}^{iq} \sum_{j=k}^{iq} \theta(j-k) \leq K^2Mq \sum_{k=0}^{q-1}  \theta (k)\,   ,
\]
and
\[
{\mathbb E}(|U_{[n/q]+1}|_{\mathbb B}^2)  \leq K^2 M \sum_{k=q[n/q] +1}^{n} \sum_{j=k}^{n} \theta(j-k)
\leq K^2M(n-q[n/q]) \sum_{k=0}^{q-1}
\theta (k)\, .
\]
Whence
\begin{equation}\label{B5}
\sum_{i=1}^{[n/q]+1}{\mathbb E}(|U_i|_{\mathbb B}^2) \leq K^2M n \sum_{k=0}^{q-1}  \theta (k)\, .
\end{equation}
Starting from \eqref{propFlo} and using  the upper bounds  \eqref{B2bis}
and \eqref{B5}, Proposition \ref{prop:maximal} follows. $\lozenge$

\subsection{A maximal inequality}

\begin{Proposition}\label{compESnmaxESn}
Let $n\geq 2$ be an integer and $Y_1, Y_2, \ldots, Y_n$ be $n$ random variables with values in a separable Banach
space $({\mathbb B}, |\cdot |_{\mathbb B})$. Assume that
 ${\mathbb P}(|Y_k|_{\mathbb B}\leq M)=1$  for
any $k \in \{1, \ldots , n\}$. Let ${\mathcal F}_1, \ldots ,
{\mathcal F}_n$ be an increasing filtration such that $Y_k$ is ${\mathcal F}_k$-measurable for any $k \in \{1, \ldots , n\}$.
Let $S_n=\sum_{k=1}^n Y_k$ and $\theta(k)$ be defined by \eqref{defthetakappendix}. Then, for any real $p>1$, the following inequality holds:
\[
{\mathbb E} \Big ( \max_{1 \leq k \leq n} |S_k|^p_{\mathbb B} \Big ) \leq  \frac 1 2 \Big ( \frac{2p}{p-1} \Big )^p \E ( |S_n |^p_{\mathbb B})+ 2^{p-1} 3^p p M^{p-1} n \sum_{k=0}^{n-2} (k+1)^{p-2} \theta( k) \, .
\]
\end{Proposition}
\noindent {\bf Proof of Proposition \ref{compESnmaxESn}.} All along the proof,  $\E_k (\cdot) = \E (\cdot | {\mathcal F}_k)$. We start by noticing that
\[
S_k = \E_k ( S_n) + \E_k ( S_k - S_n) \, .\]
Therefore
\[
\E\Big (  \max_{1 \leq k \leq n } |S_k |^p_{\mathbb B} \Big )  \leq 2^{p-1} \E \Big ( \max_{1 \leq k \leq n } |\E_k ( S_n)  |^p_{\mathbb B}  \Big ) + 2^{p-1} \E \Big ( \max_{1 \leq k \leq n } |\E_k ( S_n-S_k)  |^p_{\mathbb B}  \Big ) \, .
 \]
Notice now that $( |\E_k ( S_n)  |, {\mathcal F}_k)_{1 \leq k \leq n}$ is a submartingale. Therefore by the Doob's maximal inequality,
\[
\E \Big ( \max_{1 \leq k \leq n } |\E_k ( S_n)  |^p_{\mathbb B}  \Big )  \leq \Big ( \frac{p}{p-1} \Big )^p \E ( |S_n |^p_{\mathbb B}) \, .
\]
So, overall,
\[
 \E\Big (  \max_{1 \leq k \leq n } |S_k |^p_{\mathbb B} \Big )  \leq 2^{-1} \Big ( \frac{2p}{p-1} \Big )^p \E ( |S_n |^p_{\mathbb B})+ 2^{p-1} \E \Big ( \max_{1 \leq k \leq n } |\E_k ( S_n-S_k)  |^p_{\mathbb B}  \Big ) \, .
\]
To end the proposition, it remains to prove that
\begin{equation} \label{dec1FN}
  \E \Big ( \max_{1 \leq k \leq n } |\E_k ( S_n-S_k)  |^p_{\mathbb B}  \Big ) \leq 3^p p M^{p-1} n  \sum_{k=0}^{n-2} (k+1)^{p-2} \theta( k) \, .
  \end{equation}
With this aim, we write
\[
 \E \Big ( \max_{1 \leq k \leq n } |\E_k ( S_n-S_k)  |^p_{\mathbb B}  \Big ) = p \int_{0}^{nM}x^{p-1}  {\mathbb P} \Big ( \max_{1 \leq k \leq n } |\E_k ( S_n-S_k)  |_{\mathbb B}  > x  \Big ) dx \, .
\]
Let $q$ be a non-negative integer such that $q \leq n$. Notice that
\[
 |\E_k ( S_n-S_k)  |_{\mathbb B}  =  \Big | \sum_{i=k+1}^n \E_k (X_i)  \Big |_{\mathbb B}
 \leq  \Big | \sum_{i=k+1}^n \E_k (X_i - \E_{i-q} (X_i) ) \Big |_{\mathbb B}  +  \Big | \sum_{i=k+1}^n \E_k ( \E_{i-q} (X_i) )\Big |_{\mathbb B} \, .
\]
But
\[
 \Big | \sum_{i=k+1}^n \E_k (X_i - \E_{i-q} (X_i) )\Big |_{\mathbb B}   =  \Big |\sum_{i=k+1}^{q+k} (\E_k (X_i) - \E_{i-q} (X_i) )\Big |_{\mathbb B}   \leq 2qM \, .
\]
Therefore, for any real $x$ such that $x \in [0,n]$, choosing $q=[x]$, we get
\begin{multline*}
 {\mathbb P} \Big ( \max_{1 \leq k \leq n } |\E_k ( S_n-S_k)  |_{\mathbb B}  > 3 M x  \Big )
\leq  {\mathbb P} \Big ( \max_{1 \leq k \leq n } \Big | \sum_{i=k+1}^n \E_k ( \E_{i-[x]} (X_i) )\Big |_{\mathbb B} > M x  \Big ) \\
\leq  {\mathbb P} \Big ( \max_{1 \leq k \leq n } \E_k \Big ( \sum_{i=2}^n  |    \E_{i-[x]} (X_i) |_{\mathbb B} \Big )  > M x  \Big ) \, .
\end{multline*}
But $ (\E_k \big (\sum_{i=2}^n |\E_{i-[x]} (X_i)| \big ), \F_k )_{ 1 \leq k \leq n}$
is a martingale, so the Doob-Kolmogorov's inequality implies
\[
 {\mathbb P} \Big ( \max_{1 \leq k \leq n } \E_k \Big ( \sum_{i=2}^n  |    \E_{i-[x]} (X_i) |_{\mathbb B} \Big )  > M x  \Big )  \leq \frac{1}{Mx}
  \sum_{i=2}^n  \E \big ( |    \E_{i-[x]} (X_i) |_{\mathbb B} \big ) \leq
   \frac{n \theta( [x])}{Mx} \, .
\]
So, overall,
\begin{multline*}
 \E \big ( \max_{1 \leq k \leq n } |\E_k ( S_n-S_k)  |^p_{\mathbb B}  \big ) = p (3M)^p \int_0^{n/3} x^{p-1}{\mathbb P} \Big ( \max_{1 \leq k \leq n } |\E_k ( S_n-S_k)  |_{\mathbb B}  > 3 M x  \Big ) dx\\
\leq   3^p p M^{p-1} n  \int_0^{n/3} x^{p-2} \theta( [x]) dx \, ,
\end{multline*}
proving \eqref{dec1FN} by using the fact that $(\theta (k))_{k}$ is a non-increasing sequence. The proof of the proposition is therefore complete. $\lozenge$

\subsection{Proof of Inequality \eqref{IneThinemomentbismax}} \label{SectprInmax}

Proposition \ref{compESnmaxESn} together with Inequality
\eqref{IneThinemomentbis} leads to
\begin{multline} \label{p1inemax}
\E (\max_{1 \leq k \leq n}| S_k |^p_{\mathbb B} ) \leq 2^{-1} \Big ( \frac{2p}{p-1} \Big )^pK^p \Big ( \sum_{i=1}^n \max_{i \leq \ell \leq n } \Big \Vert | X_i |_{\mathbb B} \big | \sum_{k=i}^{\ell} \E ( X_k | {\mathcal F}_i ) \big |_{\mathbb B} \Big \Vert_{p/2} \Big )^{p/2}\\
+ 2^{p-1} 3^p p M^{p-1} n \sum_{k=0}^{n-2} (k+1)^{p-2} \theta( k) \, .
\end{multline}
Since ${\mathbb P}(|X_k|_{\mathbb B}\leq M)=1$  for
any $k \in \{1, \ldots , n\}$, it follows that
\beq \label{p2inemax}
 \sum_{i=1}^n \max_{i \leq \ell \leq n } \Big \Vert | X_i |_{\mathbb B} \big | \sum_{k=i}^{\ell} \E ( X_k | {\mathcal F}_i ) \big |_{\mathbb B} \Big \Vert_{p/2}
\leq n M^{2-2/p} \sum_{k=0}^{n-1} \theta^{2/p} (k) \, .
\eeq
On the other hand, since $(\theta(k))_{k \geq 1}$ is non-increasing,
\[
\sum_{k=1}^{n-2} (k+1)^{p-2} \theta( k) =\sum_{\ell=0}^{\log_2 (n-1)-1} \sum_{k=2^{\ell }}^{2^{\ell+1}-1}(k+1)^{p-2} \theta( k) \leq 2^{p-2}\sum_{\ell=0}^{\log_2 (n-1)} 2^{\ell (p-1)} \theta( 2^{\ell})\, .
\]
Hence, using the fact that $p\geq 2$ and again that $(\theta(k))_{k \geq 1}$ is non-increasing, we successively derive
\begin{multline*}
\sum_{k=1}^{n-2} (k+1)^{p-2} \theta( k)
\leq 2^{p-2}\Big ( \sum_{\ell=0}^{\log_2 (n-1)} 2^{\ell (2-2/p)} \theta^{2/p}( 2^{\ell}) \Big )^{p/2} \\
\leq 2^{p-2}\Big ( \theta^{2/p}(1) + 2 \sum_{\ell=1}^{\log_2 (n-1)} \sum_{k=2^{\ell -1} +1}^{2^{\ell}} 2^{\ell (1-2/p)} \theta^{2/p}( 2^{\ell}) \Big )^{p/2} \leq 2^{2p-3}\Big ( \sum_{k=1}^{n-1} k^{1-2/p} \theta^{2/p}( k) \Big )^{p/2} \, .
\end{multline*}
Since $p \geq 2$, it follows that
\beq \label{p3inemax}
\sum_{k=1}^{n-2} (k+1)^{p-2} \theta( k) \leq   2^{2p-3}n^{p/2-1}\Big ( \sum_{k=1}^{n-1} k^{1-2/p} \theta^{2/p}( k) \Big )^{p/2} \, .
\eeq
Starting from \eqref{p1inemax} and considering the upper bounds \eqref{p2inemax} and \eqref{p3inemax}, the inequality \eqref{IneThinemomentbismax} follows. $\lozenge$

\subsection{Dependence properties of Young towers}
In this section, we  assume that $T$ is a nonuniformly expanding map on $({\mathcal X}, \lambda)$ with $\lambda$ a probability measure on ${\mathcal X}$, and  that $T$ can be modelled by a Young tower. 
As in Section \ref{simple}, ${\mathcal X}$ can be any bounded  metric space and not necessarily the unit interval.
\begin{Proposition}\label{taualpha}
Let $T$ be map that can be modelled by a  Young tower with polynomial tails of the return times of order $1/\gamma$ with $\gamma \in (0,1)$. Then the inequality \eqref{majorationpourbeta}
 holds,
that is:  for any $ \alpha \in (0,1]$ there
 exists $K_\alpha>0$ such that
\begin{equation*}
   \bar \nu \Big (\sup_{f \in L_{\alpha,1}} |P^n(f) -\bar \nu (f)|
     \Big) \leq \frac{K_\alpha}{n^{(1-\gamma)/\gamma}} \, .
\end{equation*}
\end{Proposition}

\noindent {\bf Proof of Proposition \ref{taualpha}.}
 The proof is a slight modification of the proof of Theorem 2.3.6
in \cite{Go02} and is included here for the sake of completeness.
In this proof, $C$ is a positive constant,
and $C_\alpha$ is a positive constant depending only on $\alpha$. Both constants may vary
from line to line.

We keep the same notations as in Subsection \ref{YT}. For $f \in L_{\alpha}$, let
$$
  \|f\|_{L_\alpha}= L_\alpha(f) + \|f\|_\infty \, .
$$
Let $f^{(0)}= f - \bar \nu (f)$. Since $\|f^{(0)}\|_\infty \leq L_\alpha(f)$, it follows that
\begin{equation}\label{e:Vdf2}
  \|f-\bar \nu (f)\|_{L_\alpha} \leq 2 L_\alpha(f) \, .
\end{equation}

Recall that one has the decomposition
 \begin{equation}
  \label{eq_somme}
  P^n f = \sum_{a+k+b=n} \lambda_b(f) A_a ({\bf 1}_{\bar Y})  +
  \sum_{a+k+b=n} A_a E_k B_b f + C_n f \, ,
  \end{equation}
where the operators $A_n$, $B_n$, $C_n$ and $E_n$ and
 are defined in Chapter 2 of Gou\"{e}zel's PhD thesis
 \cite{Go02} and
 $\lambda_b(f)=\bar \nu(B_b(f))$.
In particular, Gou\"{e}zel has proved that
\begin{equation}\label{EandB}
\|E_k f\|_{L_\alpha}\leq \frac{C_\alpha\|f\|_{L_\alpha}} {(k+1)^{(1-\gamma)/\gamma}}  \quad \text{and} \quad
\|B_k f\|_{L_\alpha} \leq \frac{C_\alpha \|f\|_{L_\alpha}}{(k+1)^{1/\gamma}}\, .
\end{equation}

Following the proof of Lemma 2.3.5 in \cite{Go02}, there exists
a set $Z_{n}$ such that, for any bounded measurable function $g$,
\begin{equation}\label{Cn}
 |C_n(g)| \leq C\|g\|_\infty  {\bf 1}_{Z_{n}}\, ,
\end{equation}
and
\begin{equation}\label{1Cn}
\bar \nu (Z_{n}) \leq \frac{C}{ (n+1)^{(1- \gamma)/\gamma}} \, .
\end{equation}

We now turn to the term $\sum_{a+k+b=n}A_a E_k B_b f$ in
\eqref{eq_somme}. Following the proof of Lemma 2.3.3. in
\cite{Go02}, there exist a set $U_{n}$ such that, for any bounded measurable
function $g$,
\begin{equation}\label{1An}
 |A_n(g)| \leq C\|g\|_\infty  {\bf 1}_{U_{n}}\, ,
\end{equation}
and
\begin{equation}\label{1Anbis}
\bar \nu(U_{n}) \leq \frac{C}{ (n+1)^{1/\gamma}} \, .
\end{equation}
Using successively \eqref{1An} and \eqref{EandB},  we obtain that
  \begin{align}
  \label{AEB}
  \Big| \sum_{a+k+b=n} A_a E_k B_b f  \Big|  &\leq
  C \sum_{a+k+b=n}\|E_k B_b f\|_{\infty} {\bf 1}_{U_a}
  \nonumber
  \\&
 \leq  C_\alpha \sum_{a+k+b=n}\|B_b f\|_{L_\alpha}\frac {{\bf 1}_{U_a}}
  {(k+1)^{(1-\gamma)/\gamma}} \nonumber
  \\&
  \leq  C_\alpha\|f\|_{L_\alpha} \sum_{a+k+b=n}\frac {{\bf 1}_{U_a}}
  {(k+1)^{(1-\gamma)/\gamma}(b+1)^{1/\gamma}} \, .
  \end{align}

We now turn to the term
$\sum_{a+k+b=n}A_a({\bf 1}_{\bar Y})\cdot \bar \nu(B_b f)$ in
\eqref{eq_somme}. From the last equality of (2.21) in
\cite{Go02}, if $\bar \nu (f)=0$, 
\begin{align}\label{controlcentre}
  \left|\sum_{b=0}^{n-a} \bar \nu( B_b f ) \right| =
   \left| \sum_{b>n-a} \bar \nu( B_b f ) \right|
  \leq \sum_{b>n-a} \|B_b f \|_{L_\alpha}
  &\leq \sum_{b>n-a}  \frac{C_\alpha\|f\|_{L_\alpha}}{(b+1)^{1/\gamma}}\nonumber\\
  &\leq  \frac{C_\alpha\|f\|_{L_\alpha}}{(n+1-a)^{(1-\gamma)/\gamma}}\, .
  \end{align}
From \eqref{controlcentre} and  \eqref{1An}, if $\bar \nu (f)=0$,
  \begin{equation}\label{secondtermbis}
  \left|\sum_{a=0}^n A_a ({\bf 1}_{\bar Y}) \cdot \left(\sum_{b=0}^{n-a}
 \bar \nu( B_b f)\right)\right|
  \leq
  C_\alpha \|f\|_{L_\alpha} \sum_{a=0}^n \frac{{\bf 1}_{U_a}}{(n+1-a)^{(1-\gamma)/\gamma}}\, .
  \end{equation}

From \eqref{e:Vdf2}, $\|f-\nu (f)\|_{L_\alpha}\leq 2 L_\alpha(f)$. Hence, it follows from
\eqref{eq_somme},
\eqref{Cn}, \eqref{AEB} and \eqref{secondtermbis} that
\begin{equation}\label{mainbound}
|P^n(f-\nu (f))|
\leq C_\alpha L_\alpha(f) \Big ({\bf 1}_{Z_n}+ \sum_{a=0}^n \frac{{\bf 1}_{U_a}}{(n+1-a)^{(1-\gamma)/\gamma}}
+ \sum_{a+k+b=n}\frac {{\bf 1}_{U_a}}
  {(k+1)^{(1-\gamma)/\gamma}(b+1)^{1/\gamma}} \Big)\, .
\end{equation}
From  \eqref{mainbound}, \eqref{1Cn} and \eqref{1Anbis}, it follows that
\begin{multline}
\bar \nu \Big (\sup_{f \in L_{\alpha,1}} |P^n(f) -\bar \nu (f)|
     \Big)  \leq
C_\alpha \Big( \frac{1}{ (n+1)^{(1- \gamma)/\gamma}} +
\sum_{a=0}^n \frac{1}{(a+1)^{1/\gamma}(n+1-a)^{(1-\gamma)/\gamma}}
\\
 + \sum_{a+k+b=n}\frac {1}{(a+1)^{1/\gamma}
  (k+1)^{(1-\gamma)/\gamma}(b+1)^{1/\gamma}}
  \Big) \, .
\end{multline}
All the sums on right hand being of the same order (see the end of the proof of Proposition 6.2 in \cite{DDT}), it follows that there exists $K_\alpha>0$ such that
$$
\bar \nu \Big (\sup_{f \in L_{\alpha,1}} |P^n(f) -\bar \nu (f)|
     \Big) \leq \frac{K_\alpha}{n^{(1-\gamma)/\gamma}} \, ,
$$
and the proof is complete.  $\lozenge$

\subsection{Proof of Lemma \ref{constantecp}}
We shall prove here that Lemma \ref{constantecp} also holds for the derivative in the sense of Fr\'echet. Hence in the proof $D$ and $D^2$ are the first and second derivatives in the sense of Fr\'echet.

Set $|x|_q = \big ( \int_{{\mathcal X}}| x(t) |^q d \nu(t)  \big )^{1/q}$ and observe that for, any $x$ and $h$ in ${\mathbb L}^q$, by the Taylor integral formula at order 2,
\begin{multline*}
|x+h|_q^q - |x|_q^q = q \int_{{\mathcal X}}h(t) | x(t) |^{q-1} {\rm sign} ( x(t)) \mu  (dt ) \\ +
q(q-1) \int_{{\mathcal X}}h^2(t) \int_0^1 (1-s)| x(t) + s h(t) |^{q-2} ds \mu  (dt )\, .
\end{multline*}
implying that
\beq \label{1observationlem}
|x+h|_q^q - |x|_q^q = q \int_{{\mathcal X}}h(t) | x(t) |^{q-2} x(t) \mu  (dt )+ O ( |h |_q^2 ) \, .
\eeq
Define now the function $\ell$ from ${\mathbb L}^q$ to ${\mathbb R}$ by
$$
\ell (x)  = |x|_q^2 \, .
$$
Using \eqref{1observationlem}, we derive that, for any $x$ and $h$ in ${\mathbb L}^q$,
\begin{multline} \label{Dlavant}
\ell (x +h)  -\ell (x) = 2 q^{-1} \big ( \ell(x) \big )^{1-q/2} \Big (  \int_{{\mathcal X}} ( | x(t) +h(t)|^q - |x(t)|^q )  d \nu(t)\Big ) + o( ( |h |_q )\\
= 2 \big ( \ell(x) \big )^{1-q/2} \int_{{\mathcal X}}h(t) | x(t) |^{q-2} x(t) \mu  (dt )+ o(  |h |_q ) \, .
\end{multline}
Therefore $\ell$ is Fr\'echet differentiable and
\beq \label{Dl}
D \ell (x) (h) =  2 \big ( \ell(x) \big )^{1-q/2} \int_{{\mathcal X}}h(t) | x(t) |^{q-2} x(t) \mu  (dt ) \, .
\eeq
Let us prove that $\ell$ is two times Fr\'echet differentiable. Starting from \eqref{Dl}, we first write that, for any $x,h,v$ in ${\mathbb L}^q$,
\begin{multline*}
D \ell (x+v) (h) - D \ell (x) (h) = 2 \big ( \ell(x+v) \big )^{1-q/2} \int_{{\mathcal X}}h(t) | x(t) +v(t)|^{q-2} ( x(t) + v(t)) \mu  (dt )\\
 - 2 \big ( \ell(x) \big )^{1-q/2} \int_{{\mathcal X}}h(t) | x(t) |^{q-2} x(t) d\nu (t)
\end{multline*}
Notice that
\begin{multline*}
\int_{{\mathcal X}}h(t) | x(t) +v(t)|^{q-2} ( x(t) + v(t)) \mu  (dt )- \int_{{\mathcal X}}h(t) | x(t) |^{q-2} x(t) \mu  (dt )\\
= (q-1) \int_{{\mathcal X}}h(t) v(t) | x(t) |^{q-2}  \mu  (dt )+ o ( |h|_q|v|_q) \, .
\end{multline*}
Hence
\begin{multline*}
D \ell (x+v) (h) - D \ell (x) (h) = 2 (q-1) \big ( \ell(x) \big )^{1-q/2}   \int_{{\mathcal X}}h(t) v(t) | x(t) |^{q-2}  \mu  (dt ) \\
+ 2 (q-1) \big (\big ( \ell(x+v) \big )^{1-q/2} - \big ( \ell(x) \big )^{1-q/2}  \big )  \int_{{\mathcal X}}h(t) v(t) | x(t) |^{q-2}  \mu  (dt )\\
 + 2  \big ( \big ( \ell(x+v) \big )^{1-q/2}  -  \big ( \ell(x) \big )^{1-q/2} \big ) \int_{{\mathcal X}}h(t)  | x(t) |^{q-2} x(t) \mu  (dt )+ o ( |h|_q|v|_q) \, .
\end{multline*}
Using \eqref{Dlavant}, we infer that
$$
(\big ( \ell(x+v) \big )^{1-q/2} - \big ( \ell(x) \big )^{1-q/2}  = (2-q) \big ( \ell(x) \big )^{1-q}  \int_{{\mathcal X}}v(t) | x(t) |^{q-2} x(t) \mu  (dt )+ o(|v|_q) \, .
$$
So, overall,
\begin{multline*}
D \ell (x+v) (h) - D \ell (x) (h) = 2 (q-1) \big ( \ell(x) \big )^{1-q/2}   \int_{{\mathcal X}}h(t) v(t) | x(t) |^{q-2}  \mu  (dt ) \\
 + 2 (2-q) \big ( \ell(x) \big )^{1-q}  \int_{{\mathcal X}}v(t) x(t) | x(t) |^{q-2}  \mu  (dt ) \int_{{\mathcal X}}h(t) x(t) | x(t) |^{q-2} \mu  (dt )+ o ( |h|_q|v|_q) \, .
\end{multline*}
Therefore $\ell$ is two-times Fr\'echet differentiable and
\begin{multline} \label{Dl2}
D^2\ell (x) (h,v) = 2 (q-1) \big ( \ell(x) \big )^{1-q/2}   \int_{{\mathcal X}}h(t) v(t) | x(t) |^{q-2}  \mu  (dt ) \\
 + 2 (2-q) \big ( \ell(x) \big )^{1-q}  \int_{{\mathcal X}}v(t) x(t) | x(t) |^{q-2}  \mu  (dt ) \int_{{\mathcal X}}h(t) x(t) | x(t) |^{q-2} \mu  (dt ) \, .
\end{multline}
Since $\psi_p (x) = \big ( \ell (x) \big )^{p/2}$, $\psi_p$ is also two-times Fr\'echet differentiable. Moreover
\[
D \psi_p (x) (h) =  2 \big ( \ell(x) \big )^{(p-q)/2} \int_{{\mathcal X}}h(t) | x(t) |^{q-2} x(t) \mu  (dt ) \, ,\]
and
\begin{multline} \label{Dl2psip}
D^2 \psi_p (x) (h,v) = p (q-1) \big ( \ell(x) \big )^{(p-q)/2}   \int_{{\mathcal X}}h(t) v(t) | x(t) |^{q-2}  \mu  (dt ) \\
 + p (p-q) \big ( \ell(x) \big )^{-q+p/2}  \int_{{\mathcal X}}v(t) x(t) | x(t) |^{q-2}  \mu  (dt ) \int_{{\mathcal X}}h(t) x(t) | x(t) |^{q-2} \mu  (dt ) \, .
\end{multline}
Starting from \eqref{Dl2psip} and using the fact that $\ell (x) = |x|_q^2$, we get
\begin{multline} \label{Dl2psipbis}
D^2 \psi_p (x) (h,v) = p (q-1) | x |_q^{p-q}  \int_{{\mathcal X}}h(t) v(t) | x(t) |^{q-2}  \mu  (dt ) \\
 + p (p-q) |x|_q^{p-2q}  \int_{{\mathcal X}}v(t) x(t) | x(t) |^{q-2} \mu  (dt ) \int_{{\mathcal X}}h(t) x(t) | x(t) |^{q-2} \mu  (dt ) \, ,
\end{multline}
and an application of H\"older's inequality shows that ${\mathbb L}^q $ belongs  to the class ${\widetilde {\mathcal C}}_2 (p, {\tilde c}_p)$ with ${\tilde c}_p = p \big ( \max( p, 2q-p) -1\big )$. To prove that ${\mathbb L}^q $ belongs  to the class
${\ {\mathcal C}}_2 (p, c_p)$ with $c_p=p(\max(p,q)-1)$, it suffices to write \eqref{Dl2psipbis}  with $h=v$, and to use the fact that $\big ( \int_{{\mathcal X}}v(t) | x(t) |^{q-2} x(t) \mu  (dt )\big )^2$ is non-negative. This ends the proof of Item 1.

The proof of Item 2  is omitted since it uses the same arguments as for ${\mathbb L}^2$. $\lozenge$

\end{document}